\documentclass[journal]{IEEEtran}

\IEEEoverridecommandlockouts 

\usepackage{graphicx}
\usepackage{amsthm}
\usepackage{amssymb,amsmath,amsfonts}
\usepackage{color}
\usepackage[dvipsnames]{xcolor}
\usepackage{epsfig}
\usepackage{cases}
\usepackage{epstopdf}
\usepackage{cases}
\usepackage{cite}
\usepackage[colorlinks,
linkcolor=blue,
anchorcolor=green,
citecolor=blue
]{hyperref}

\newtheorem{theorem}{Theorem} 

\newtheorem{lemma}{Lemma}

\newtheorem{definition}{Definition}

\newtheorem{remark}{Remark}

\newtheorem{assumption}{Assumption}

\newcommand{\nm}{\nonumber}

\newcommand{\bR} { {\mathbb R}}
\newcommand{\bN} { {\mathbb N}}
\newcommand{\1}{\mbox{1}\hspace{-0.25em}\mbox{l}}

\begin{document}
	\title{Fixed-Time Cooperative Behavioral Control for Networked Autonomous Agents with Second-Order
		Nonlinear Dynamics} 
	\author{Ning~Zhou, Xiaodong~Cheng, Zhongqi Sun, Yuanqing~Xia
		\thanks{*This work was supported in part by the National Natural Science Foundation of China under Grant 61603095, Grant 61972093 and Grant 61720106010. The work of Yuanqing Xia was also supported in part by the Science and Technology on Space Intelligent Control Laboratory under Grant KGJZDSYS-2018-05. \textit{(Corresponding author: Xiaodong Cheng.)}}% <-this % stops a space
		\thanks{Ning Zhou is with {School of Electrical Engineering, Hebei University of Science and Technology, Shijiazhuang 050018, China.}
			{\tt\small zhouning2010@gmail.com}.}%%
		\thanks{Xiaodong Cheng is with the Department of Engineering, University of Cambridge, Trumpington Street, Cambridge, CB2 1PZ, United Kingdom
			{\tt\small  xc336@cam.ac.uk}.}%
		\thanks{Zhongqi Sun and Yuanqing Xia are with the School of Automation,
			Beijing Institute of Technology, Beijing 100081, China.
			{\tt\small zhongqisun@bit.edu.cn, xia\_yuanqing@bit.edu.cn}.}
	}
	
	\maketitle

	\begin{abstract}
		In this paper, we investigate the {fixed-time} behavioral control {problem for a team} of {second-order nonlinear agents}, {aiming} to achieve \textit{{a desired} formation} {with} \textit{collision/obstacle~avoidance}. {In the proposed approach, the two behaviors(tasks) for each agent are prioritized and integrated via the framework of the null-space-based behavioral projection, leading to a desired merged velocity that guarantees the fixed-time convergence of task errors. To track this desired velocity, we design a fixed-time sliding mode controller for each agent with state-independent adaptive gains, which provides a fixed-time convergence of the tracking error. The control scheme is implemented in a distributed manner, where each agent only acquires information from its neighbors in the network. Moreover, we adopt an online learning algorithm to improve the {robustness} of the closed system with respect to uncertainties/disturbances. Finally, simulation results are provided to show the effectiveness of the proposed approach.}
	\end{abstract}
	
	
	\section{Introduction}
	Recently, intelligent multi-robot systems  have found broad applications in e.g., cooperative reconnoiter, monitoring and rescue missions {\cite{sargolzaei2020control}}, formation of autonomous robots\cite{oh2015survey}, coordination of spacecraft {\cite{zou2020distributed}}, and {coordinated path-following of surface vessels \cite{Zhang2022Constrained}}. These systems that bond multiple autonomous agents (e.g., vehicles and robots) by communication networks can carry out much more complicated tasks than those that a single agent can ever accomplish. 
	However, in many real applications, autonomous agents are deployed in a complex and dynamic environment to execute multiple parallel tasks, e.g., to maintain a desired formation and avoid moving obstacles simultaneously. How to operate such systems efficiently and safely poses a challenging control problem.
	
	To resolve multi-mission control problems for multi-agent systems, the so-called behavioral approach is developed, see e.g., \cite{brooks1986robust,antonelli2008null} and the reference therein.
	In this scheme, {a comprehensive motion task is decomposed into multiple smaller and simpler subtasks,  described by \textit{behavioral functions}, which {generate} a set of motion commands. The eventual motion control of each individual agent is performed as the outcome of merging multiple behaviors simultaneously. To  merge multiple prioritized subtasks, low-priority tasks are projected to the null space of higher-priority tasks.   
		This scheme is also referred to \textit{null-space-based behavioral approach} (NSB), and in the works \cite{antonelli2006kinematic,antonelli2008stability,antonelli2008null,antonelli2010flocking, schlanbusch2011spacecraft,huang2019adaptive,zhou2018neural}, different centralized NSB approaches are developed for controlling a team of autonomous vehicles to cooperatively carry out multiple tasks.
		However, all the above-mentioned behavioral control approaches are formulated in a centralized fashion, which requires the global information of overall multi-agent systems. In many real-world applications, information acquisition in a global level may not be practical due to communication cost and constraints.} 
	
	In those scenarios, distributed control schemes are resorted, in which controllers are localized at each autonomous agents in a network. Through information exchange among neighboring agents, the networked agents can accomplish certain cooperative tasks together. This concept of distributed control has shown a great potential in various applications of multi-agent systems, see e.g., \cite{ChengACOM2018Power,zou2019event,liu2019event} for an overview. 
	{To achieve formation control of networked robots in an environment with obstacles, different methods have been developed, including artificial potential field \cite{nguyen2016formation}, geometric optimization \cite{alonso2016distributed}, fluid-based approach \cite{wu2019formation}, and  model predictive control methods \cite{nascimento2016multi,dai2017distributed}.  However, the control configuration in these methods {is} restricted to perform only two subtasks, namely, obstacle avoidance and formation, at the same time. To handle multiple tasks for a team of robots, behavioral approaches that are implemented in a distributed/decentralized way present a promising direction. A decentralized framework of behavioral approaches was firstly given in \cite{antonelli2010flocking}, although {a} theoretic guarantee of the convergence of behavior errors is lacking. A distributed  formation control method using NSB is provided in \cite{ahmad2014behavioral}, which results in the asymptotic stability of the closed-loop system. However, this method is limited to triangular formation in {an} obstacle-free environment. 
		
		In contrast to the existing literature, this paper presents a new control framework that combines the concept of \textit{fixed-time control} with behavioral approaches, and we also provide a distributed implementation of this framework. The {benchmarking} work on fixed-time control for generic nonlinear systems was presented in \cite{polyakov2012nonlinear}, and it shows a fast convergence rate, high-precision control performance, and disturbance rejection properties \cite{zuo2018overview}. These merits are inherited by our fixed-time behavioral based framework. Moreover, in contrast to our previous finite-time methods in \cite{zhou2017finite}, the convergence time (settling time) of our procedure can be predicted without requiring knowledge of the initial conditions. 
		
		In this work, we apply the proposed framework to control networked multi-agent systems to achieve multiple tasks in a fixed time. The major challenges in this framework are to provide a theoretical guarantee on the fixed time property for multiple tasks, particularly when tasks are conflicting with each other, and to handle the {local minima} when the desired velocities for different tasks cancel out each other. Moreover, we extend the preliminary results of this paper in \cite{Zhou2020FTB} by considering rather general settings for the agent dynamics, which is modeled by a second-order dynamics including nonlinear uncertainty and unknown external disturbances. In a dynamic environment with moving obstacles, the networked agents need to achieve a certain formation while avoiding collisions with each other and the obstacles.} To the best of our knowledge, solving such a problem in a fixed-time setting has not been addressed by any existing methods so far. To solve the problem, we introduce the behavior functions of collision avoidance and cooperative formation, respectively, which are merged in priority via the null-space-based behavioral projection to give the desired velocity for each agent that can be computed based on only local information. We then design distributed fixed-time controllers for the agents to cooperatively track the fixed-time desired velocities, where the universal approximation property of the radial basis function neural networks (RBFNNs) is applied to identify the uncertain terms in the system. Specifically, the contributions
	of this paper are emphasized as follows:
	
	(1) The concept of a fixed-time control scheme is used for the first time in the framework of behavioral approach. Based on a new designed distributed fixed-time estimator, a distributed fixed-time behavioral strategy is designed at the kinematics level, which leads to cooperative behaviors of multi-agent systems with second-order nonlinear dynamics. The developed fixed-time behavior approach can handle various shapes of both flexible and {fixed} formations in a distributed framework and guarantee collision/obstacle avoidance.
	
	(2)  A kind of state-independent adaptive gain is designed and contributed to constructing a set of fixed-time intelligent tracking control laws, {which allows for an adjustable control accuracy even after the settling time. Moreover, we provide a theoretical guarantee on the fixed-time convergence of velocity and position tracking errors} and strengthen the robustness of the closed-loop system.
	
	{The rest of the paper is organized as follows. In Section~\ref{sec2}, we provide some preliminaries on the fixed-time control and behavioral approach. Then, the control problem is formulated; In Section~\ref{sec3}, the desired velocity for each agent is designed using the fixed-time behavioral control scheme, and Section~\ref{sec4} provides the controller to track the desired velocity; {Two adjustable  control gains in the controller are discussed in Section~\ref{rem5}}; Section~\ref{sec5} shows the simulation results in three-dimensional space. Finally, Section~\ref{sec:Conclusions} concludes the paper.}
	
	\textit{Notation:} The set of real numbers is denoted by $\bR$. For a vector or matrix, $\| \cdot \|$ denotes its Euclidean norm. The $i$-th element of a vector $v$ is denoted by $v_{i}$. The operator $ \mathrm{blkdiag}\{\cdot\} $ defines a block diagonal matrix. $ x^{[p]}:=|x|^p\mathrm{sgn}(x), x, p\in \bR $. $\mathrm{sgn}(\cdot )$ is the sign function that returns $-1$, $0$ or $1$.  
	
	\section{Preliminaries and Problem Formulation}\label{sec2}
	This section presents the preliminaries with regard to the fixed-time stability and prioritized multi-behavior composition. Then, the problem is formulated for the cooperative control of networked   agents with second-order nonlinear dynamics.

	\subsection{Fixed-Time Stability}
	
	Consider a nonlinear system
	\begin{align}\label{system0}
	\dot{x}(t)=f(t,x),\ x(0)=x_{0},
	\end{align}
	with $ x(t) \in \bR^{n} $ and the nonlinear function $f(t,x)$. If $ f(t,x) $ is discontinuous, the solutions of \eqref{system0} are Filippov. Suppose the origin is an equilibrium point of \eqref{system0}, then the fixed-time stability is defined as follows.
	
	\begin{definition} \cite{polyakov2012nonlinear}
		The origin $x = 0$ is said to be \textbf{globally fixed-time stable} if it is globally asymptotically stable and any solution $x(t,x_0)$ of \eqref{system0} reaches $x = 0$ in some settling time $t = T(x_0)$ and remains there for all $t \geq T(x_0)$, where $T(x_0)$ is globally bounded by some number $T_{\max} \in \mathbb{R}_{>0}$. 
	\end{definition}
	
	Notice that in the concept of the fixed-time stability, the settling (convergence) time $T(x_0)$ is always bounded independent of the initial condition $x_0$. In the terms of the Lyapunov stability theory, the fixed-time stability of the nonlinear system \eqref{system0} can be characterized by the following lemma.
	\begin{lemma}\cite{polyakov2012nonlinear,parsegov2013fixed}
		\label{ld1} 
		If there exists a continuous radially unbounded and positive definite function $ V:\bR^{n}\to\bR_{>0}$ such that $ V(x)=0$ if and only if $x=0 $, and any solution $ x(t,x_0) $ of \eqref{system0} satisfies
		\begin{align}\label{V0}
		\dot{V}(x)&\le -\eta_{1}V^{k_{1}}(x)-\eta_{2}V^{k_{2}}(x),\\
		{\rm or}\ \dot{V}(x)&\le -(\eta_{1}V^{k_{3}}(x)+\eta_{2}V^{k_{4}}(x))^{k_{5}},\label{V01}
		\end{align}
		where $ \eta_{1},\eta_{2},k_{1},k_{2},k_{3},k_{4},k_{5} \in \mathbb{R}_{>0}$ with $ k_{1}>1 $, $ 0<k_{2}<1 $, $ k_{3}k_{5}>1 $, and $ k_{4}k_{5}<1 $, then the origin of \eqref{system0} is globally fixed-time stable and the settling time function $ T $ can be estimated by 
		\begin{align*}
		&T\le T_{\max}:=\frac{1}{\eta_{1}(k_{1}-1)}+\frac{1}{\eta_{2}(1-k_{2})}, \\
		\text{or} \ & T\le T_{\max}:= \frac{1}{\eta_{1}^{k_{5}}(k_{3}k_{5}-1)}+\frac{1}{\eta_{2}^{k_{5}}(1-k_{4}k_{5})},
		\end{align*}
		where $T_{\max}$ is independent on the initial condition $x(0)$. 
	\end{lemma}
	
	\begin{remark}\label{remark2}
		Let $\eta_{0} \in \mathbb{R}_{>0}$. If we replace \eqref{V0} and \eqref{V01} in Lemma~\ref{ld1} by
		\begin{align*}
		&\dot{V}(x)\le -\eta_{0}V(x)-\eta_{1}V^{k_{1}}(x)-\eta_{2}V^{k_{2}}(x), 
		\\
		{\text{and}}\ & \dot{V}(x)\le-\eta_{0}V(x) -(\eta_{1}V^{k_{3}}(x)+\eta_{2}V^{k_{4}}(x))^{k_{5}},
		\end{align*}
		respectively, then the conclusion of Lemma~\ref{ld1} still hold, due to the fact that $ -\eta_{0}V(x)\le0 $. 
	\end{remark}
	{\subsection{Prioritized Multi-Behavior Composition}
		
		Generally, a behavior (mission/task) involving some agents may require the simultaneous accomplishment of several submissions.
		The NSBC uses a geometric hierarchical composition of the behaviors' outputs to obtain
		motion-reference signals for each agent \cite{baizid2017behavioral,antonelli2006kinematic}. The scheme has three levels:
		1) \emph{Elementary behaviors/missions} {are} the fundamental mission/task functions to be controlled in the kinematic level;
		2) \emph{Composite behaviors/missions} are the combinations of elementary behaviors in a prioritized order;
		3) \emph{Supervisor} is used to switch between the defined composite behaviors/missions.

		Let $\rho_k: \bR_{\ge 0}\to\mathbb{R}^{m_k}$ be the behavior function for any $\ 1 \leq k \leq r $, {where $k\in \bN$ denotes the $ k $th behavior, $r \in \bN$ is the total number of the behaviors, and $m_k\in \bN_{\ge1}$ is the dimension of $ k $th behavior space}. Then we define a behavior hierarchy which complies with the following rules:
		\begin{enumerate}
			\item Assume that $k=1$ is the top priority. Here $k_a < k_b$ means that $k_a$ is higher in priority than $k_b$. A behavior of priority $k_b$ may not disturb the other behavior of priority $k_a$. The lower priority behaviors are executed in the null space of all higher priority behaviors.
			
			\item For any $ \ 1 \leq k \leq r $, the behavior Jacobian matrix $J_k\in \mathbb{R}^{m_k \times mn}$ determines the mappings from the joint velocities to the behavior velocities, where $ m $ is the dimension of all the system states, $n $ denotes the number of agents.
			\item The dimension of the lowest level behavior $m_r$ may be larger than $mn-\sum_{k=1}^{r-1}m_k$ so that the dimension $mn$ of the joint space exceeds the entire dimension of all behaviors.
	\end{enumerate}}
	\subsection{Problem Formulation}
	We consider a group of {$n\ (n \geq 2)$} {autonomous agents with second-order nonlinear dynamics described as }
	\begin{subequations} \label{agent}
		\begin{align}
		\dot{x}_{i1}(t) =& x_{i2}(t) \label{1a},\\
		\dot{x}_{i2}(t) =& u_{i}(t)+{f_{i}(\bar{{x}}_{i})+d_{i}(\bar{{x}}_{i},t)},\ i = 1, \ldots, n,\label{1b}
		\end{align}
	\end{subequations}
	where $x_{i1}, x_{i2}: \bR_{\ge 0} \to \bR^{3}$ are the position and velocity state vectors of the $i$th agent, respectively. {The stacked vector $\bar{x}_{i}:=[x_{i1}^{\top}, x_{i2}^{\top}]^{\top}\in \mathbb{R}^{6}$ then represents the overall states of the $ i $th agent.}
	$u_{i}: \bR_{\ge 0} \to \bR^{3}$ is the control input, and ${f_{i}({\bar{x}}_{i})}: \bR^{6} \to \bR^{3}$ {is the unknown uncertainty, which is} locally Lipschitz with $f_{i}(0)=0$. $d_{i}(\bar{{x}}_{i},t): \bR^{6} \to \bR^{3}$ is an unknown external disturbance to the system \eqref{agent}. 
	
	{The $n$ agents are coupled via a communication network described by an undirected weighted graph $G$ with the node set ${V}: = \{1,2, \cdots, n\}$. Let $A = [a_{ij}]\in\bR^{n\times n}$ be the weighted adjacency matrix of $G$, and $ a_{ij} $ where $a_{ij}=a_{ij}\in\bR_{\ge0}$  denotes the communication strength between agents $i$ and $j$. {Denote $ D:={\rm diag} \{d_{11},\dots,d_{nn}\} $ and $ L\in\bR^{n\times n}$ be the degree matrix and Laplacian matrix of graph $ G $, where $ d_{ii}=\sum_{j=1}^{n}a_{ij} $ for $ i=1,\dots,n $ and $ L=D-A $. Consider an auxiliary graph $ G' $, which represents the interactions among a virtual leader and $ n $ agents as followers. The \textit{leader adjacency matrix} is defined by $ B={\rm{diag}}\{b_{1},\dots,b_{n}\}\in\bR^{n\times n} $, where $ b_{i} > 0$ if the information of the virtual leader is available to the follower agent $i$, and $ b_{i} = 0 $ otherwise. 
			\begin{assumption}\label{a2.1}
				Assume that at each time instant, there exists at least an agent $i \in V$ such that $ b_{i} =1 $.
	\end{assumption}}}
	
	{Consider an unknown environment with dynamic obstacles. This paper aims to design a distributed behavioral control scheme for the networked agents to form a predefined formation within a fixed time and meanwhile avoid colliding with each other as well as environmental obstacles. 
	}
	\section{Fixed-Time Behavioral Control Design}\label{sec3}
	In this section, two types of behaviors are analyzed, namely, the collision-avoidance behavior and cooperative behavior, {which yield two desired velocities. Both velocities guarantee fixed-time convergence, and they are then prioritized and merged to a desired velocity} for the follow-up tracking control.
	\subsection{Collision Avoidance Behavior}\label{sec3.1}
	This paper considers a dynamic environment, which allows for environmental obstacles with time-varying positions. It is required for each agent to avoid both dynamic environmental obstacles and the other moving agents. This is referred to as the Collision Avoidance Behavior {(CoAB)}, {which can be characterized by a {CoAB} function for each individual agent. Particularly,  the {CoAB} function is defined as the shortest distance, in terms of time, between an agent and all the other objects, including both the environmental obstacles and the other agents}. 
	When a team of autonomous agents
	encounters an obstacle, {this {CoAB} function can be used to generate a desired} driving velocity to keep each agent maneuvering at a safe
	distance from all the other objects. 
	
	Let $x_{i}^o:\bR_{>0}\to\bR^{3}$ be the position of the closest object to the agent $i$, that yields the {CoAB} function $\rho_{io}:\bR^{3}\to\bR_{>0}$ as 
	\begin{equation}\label{1}
	\rho_{io} = \frac{1}{2}\|x_{i1}-x_{i}^{o}\|^{2}.
	\end{equation}
	Then the behavior-dependent Jacobian matrixes are given as
	\begin{align}\label{jio}
	J_{io} = & \frac{\partial \rho_{io}}{\partial x_{i1}} ={(x_{i1}-x_{i}^o)^{\top}} \in \mathbb{R}^{1 \times 3},\\  J_{i}^{o} = &\frac{\partial \rho_{io}}{\partial x_{i}^{o}} =-{(x_{i1}-x_{i}^o)^{\top}} \in \mathbb{R}^{1 \times 3},\nm
	\end{align}
	with the right pseudoinverse
	$ J_{io}^{\dag} =({x_{i1}-x_{i}^o})/{\|(x_{i1}-x_{i}^o)\|^{2}}, 
	$ that represents the unit vector aligned along the direction from the nearest object to the agent $i$. We define a circular \textit{repulsive zone} around each object, with the coordinates of the object as the center and $d\in\bR_{>0}$ as the radius. Then, we define the {CoAB} task error as 
	\begin{equation}\label{2}
	\tilde{\rho}_{io} := \rho_{od}-\rho_{io}, \ \text{with}~\rho_{od}: = \frac{d^{2}}{2},
	\end{equation}%
	from which the desired collision-avoidance behavioral velocity $\dot{x}_{io} :\bR_{\ge0}\to\bR^{3} $ is designed for the agent $ i $ as
	\begin{align}\label{v1}
	\dot{x}_{io} =& J_{io}^{\dag} [\lambda_{io}\alpha_{io}(\tilde{\rho}_{io})-J_{i}^{o}\dot{x}_{i}^{o}].
	\end{align}
	In \eqref{v1}, $\lambda_{io} \in \bR_{>0}$ is a state-dependent gain to be determined, {$ \dot{x}_{i}^{o} $ is the velocity of the obstacle}, and $\alpha_{io}(\tilde{\rho}_{io}):\bR_{>0}\to \bR^{3}$ is a {continuous} and differentiable function designed as follows:
	\begin{align} 
	&\alpha_{io}(\tilde{\rho}_{io}):= \label{6}\\ \nm
	&\begin{cases}
	\wp_{1}\tilde{\rho}_{io}+\wp_{2}\tilde{\rho}_{io}^{[2]}, &\mathrm{if}\ \sigma_{1,i} (\dot{\tilde{\rho}}_{io},\tilde{\rho}_{io})\neq0,\ |\tilde{\rho}_{io}| \leq \phi_{s},\\
	(\beta_{1}\tilde{\rho}_{io}^{[r_{1}]}+\beta_{2}\tilde{\rho}_{io}^{[r_{2}]})^{[r_{0}]},&\mathrm{otherwise},
	\end{cases}
	\end{align}
	with
	\begin{align} 
	&\sigma_{1,i} (\dot{\tilde{\rho}}_{io},\tilde{\rho}_{io}):=\dot{\tilde{\rho}}_{io}+c_{0}(\beta_{1}\tilde{\rho}_{io}^{[r_{1}]}+\beta_{2}\tilde{\rho}_{io}^{[r_{2}]})^{[r_{0}]}, \nonumber\\
	&\wp_{1}:=(2-r_{0})
	(\beta_{1}\phi_{s}^{[\frac{r_{1}r_{0}-1}{r_{0}}]}+\beta_{2}\phi_{s}^{[\frac{r_{2}r_{0}-1}{r_{0}}]})^{[r_{0}]},\label{8} \\ &\wp_{2}:=(r_{0}-1)(\beta_{1}\phi_{s}^{[\frac{r_{1}r_{0}-2}{r_{0}}]}+\beta_{2}\phi_{s}^{[\frac{r_{2}r_{0}-2}{r_{0}}]})^{[r_{0}]},\label{9}
	\end{align}
	where $\beta_{1},\beta_{2}, c_0, \phi_{s}, r_{0},r_{1},r_{2} \in \bR_{>0} $ are design parameters with $r_{1}r_{0}>1$, $r_{2}r_{0}<1$, and $r_{0}\in(\frac{1}{2},1) $. The form of the function \eqref{6} is originally used to design the fixed-time terminal sliding mode (FTTSM) in \cite{jiang2016fixed}, while this paper adapts this FTTSM function with an additional {parameter $c_{0}$, which effects the convergence time of task error in terms of fixed-time control theory}. 	{In \eqref{6}, the larger values of $\lambda_{io}$, $r_1$, $\beta _1$ and $\beta _2$, or the smaller values of $r_2$ and $\phi _s$, lead to a faster convergence rate of the task errors. However, they may also lower the convergence accuracy of the robust tracking control, see Remark~\ref{rem:parameter}.}
	{Moreover, the coefficients $ \wp_{1} $ and $ \wp_{2} $ in the function $ \alpha_{io}(\tilde{\rho}_{io}) $ are designed to make $ d\alpha_{io}(\tilde{\rho}_{io})/dt $ as a continuous function of the time for $ i=1,\dots,n $, as well as functions \eqref{if} and \eqref{32} in the following design.}

	\begin{remark}
		In \cite{schlanbusch2011spacecraft}, the desired velocity is designed as $\dot{x}_{io} = J_{io}^{\dag} \lambda_{io}\tilde{\rho}_{io}$, with $\lambda_{io}$ a positive scalar gain. However, such a design only addresses static obstacle and provides asymptotic convergence of the task error. In contrast, with the desired velocity $ \dot{x}_{io} $ in \eqref{v1}, we can handle moving obstacles and ensure a fixed-time convergence of the task error. {Furthermore, the proposed collision avoidance scheme \eqref{v1} requires the information of obstacles' position and velocity. In practical application, the positions of obstacles can be obtained via  lidar, camera or ultrasonic sensors, etc., and then their velocities can be calculated (or estimated) by taking the time derivative of their positions.}
	\end{remark}
	
	{In the following lemma, we show that the designed desired velocity 
		\eqref{v1} guarantees the collision avoidance for each agent. }%
	\begin{lemma}\label{l1}
		Consider the collision avoidance behavior function \eqref{1} for each agent $i \in V$. Suppose that the agent $i$ is driven by the desired velocity \eqref{v1}. If $ \|x_{i1}-x_{i}^{o}\|\le d $, then there exists a globally bounded settling time $T_{i,o} \in \mathbb{R}_{>0}$ such that $\|x_{i1}-x_{i}^{o}\|\ge \tilde{d} $, for any initial task error $\tilde{\rho}_{io}(0)\in\bR$, $0 < \phi_{s} \leq d^2/2$ and for all $t\geq T_{i,o}$, {where $\tilde{d}:=\sqrt{d^{2}-2\phi_{s}}$}.
	\end{lemma}
	Note that by properly {choosing} the parameter $\phi_s$, the value $\tilde{d}$ can be set as the safe distance from each agent to its nearest object. Then the designed velocity \eqref{v1} guarantees that if the agent $i$ enters the repulsive zone, it will be kept away from the object with the safe distance. If the agent $i$ is not inside the repulsive zone of any other object, i.e., $ \|x_{i1}-x_{i}^{o}\|> d $, then the collision avoidance task is not active.
	
	\subsection{Cooperative Behavior with Fixed-Time Estimator}\label{sec3.2}
	This section considers the cooperative behavior, where all the agents are moving towards a predefined formation. To achieve this task, we resort to the virtual leader approach \cite{egerstedt2001virtual,porfiri2007tracking}. In this scheme, a virtual leader is specified as
	a reference for the networked agents, which are designed to maintain a desired offset with respect to the position of the
	virtual leader. Thereby,
	we define the Cooperative Tracking Behavior (CTB) function $\rho_f: \bR^{3n}\to \bR^{n}$ as
	\begin{align} \label{ctbf}
	\rho_f = [\rho_{1f},\dots,\rho_{nf}]^{\top}, \ \text{with} \
	\rho_{if} := \frac{1}{2}\|x_{i1}-\hat{x}_{i1}\|^{2},
	\end{align}
	where $\hat{x}_{i1} \in \mathbb{R}^3$ denotes a pre-estimated position for the agent $i$, which is supposed to be driven towards $\hat{x}_{i1}$ and kept a distance $d_{i0}\in\bR_{>0}$ from $\hat{x}_{i1}$.
	Then, we denote the cooperative tracking behavior error as
	\begin{align}\label{15}
	{\tilde{\rho}_{f}:=[\tilde{\rho}_{1f},\dots,\tilde{\rho}_{nf}]^{\top}, \ \text{with} \ \tilde{\rho}_{1f}:=\frac{d^2_{i0}}{2}-{\rho}_{if}}.
	\end{align} 
	It is worth to emphasize that the design of the pre-estimated position $\hat{x}_{i1}$ is essential for the coordination accuracy and can lead to a distributed cooperative control scheme. 
	
	{First, we provide a distributed {fixed-time} observer for each agent to estimate the position of the virtual leader. Denote}
	\begin{equation} \label{eq:bar_x}
	\bar{x}_{i1}(t)=\hat{x}_{i1}(t)-x_{o}(t), \ {\bar{x}_{01}=0, \ i=1,\dots,n,} 
	\end{equation}
	{with $x_o:\bR_{>0}\to\bR^{3}$, a second order differentiable function, the real position of the virtual leader.}
	Then, we design the trajectory for $\hat{x}_{i1}$ in \eqref{11} using a fixed-time sliding mode estimator as follows:
	\begin{align}
	\dot{\hat{x}}_{i1}(t)=& -K_{1}\eta_{i1}(t)-K_{2}\eta_{i2}(t)- K_{3}{\eta_{i3}(t)},\label{11}\\
	{\eta}_{i1}(t)=& {({\sum_{j=0}^{n}}a_{ij}(\bar{x}_{i1}(t)-\bar{x}_{j1}(t)))^{[\frac{r_{3}}{r_{4}}]},}\nm\\
	{\eta}_{i2}(t)=& {(\sum_{j=0}^{n}a_{ij}(\bar{x}_{i1}(t)-\bar{x}_{j1}(t)))^{[\frac{r_{5}}{r_{6}}]},}\nm\\
	{\eta}_{i3}(t)=& {{\rm{sgn}}(\sum_{j=0}^{n}a_{ij}(\bar{x}_{i1}(t)-\bar{x}_{j1}(t))),}\nm
	\end{align}
	where {$a_{ij}$ is the entry of the adjacency matrix $A$ indicating the coupling strength between the agents $ i $ and $ j $, and $ a_{i0}=b_{i} $.} $K_{1},K_{2},K_3\in\bR_{>0}$ are estimator gains to be designed, {$r_{3},r_{4},r_{5},r_{6}\in\bR_{>0} $ are design parameters satisfying $ \frac{r_{3}}{r_{4}}>1 $ and $\frac{ r_{5}}{r_{6}}<1 $. {Moreover, the larger values of $ K_{1} $, $ K_{2} $, $ r_{3} $, and $ r_{6} $, or smaller values of $ r_{4} $ and $ r_{5} $, give a faster convergence speed of the estimation error but a lower estimation precision.}

		With the design of the estimator $\hat{x}_{i1}$ in \eqref{11} for each agent $i$, the following lemma claims that \eqref{11} is a fixed-time estimator of the leader's states $ x_{o} $, namely after a fixed time, $\hat{x}_{i1}(t)$ converges to $ x_{o}(t)$. 
		\begin{lemma}\label{l8}
			Consider the estimator \eqref{11} with the gains $K_1, K_2, K_3\in\bR_{>0} $. Define $H: = L + B$, where $L$ and $B$ are the Laplacian matrix of $G$ and the leader adjacency matrix of $G'$, respectively. If the velocity of the  leader is bounded as 
			\begin{equation}\label{a3.1}
			{K_{3}\ge  \sup_{t\ge 0}\|\dot{x}_{o}(t)\|_{\infty},}
			\end{equation}
			{then there exists} a globally bounded settling time $T_e \in \mathbb{R}_{>0}$
			such that $\hat{x}_{i1}(t)\equiv x_{o}(t) $ for any initial condition $ (\hat{x}_{i1}(0),x_{o}(0))\in \bR^{3} \times \bR^{3} $ and for all $ t\ge T_{e} $. 
	\end{lemma}}

	Next, we use the fixed-time estimator in \eqref{11} to design the desired velocity of each agent $i$ for {eliminating} the cooperative behavior error in \eqref{15}. The design procedure follows similarly as in Section~\ref{sec3.1}. First, we calculate the Jacobian matrixes $ J_{f} $, $ J_{\hat{f}} $ and the right pseudo-inverse of $ J_{f} $ based on the CTB function \eqref{ctbf} as follows:
	\begin{align} \label{jif}
	{J_{f}}=&{\rm blkdiag}\left\{\frac{\partial \rho_{if}}{\partial x_{11}},\dots,\frac{\partial \rho_{if}}{\partial x_{n1}}\right\}
	\\ 
	\nm=&{\rm blkdiag}\left\{{{(x_{11}-\hat{x}_{11})^{\top}},\ldots, {(x_{n1}-\hat{x}_{n1})^{\top}}}\right\}\in\mathbb{R}^{n\times
		3n},
	\\
	J_{\hat{f}}=&{\rm blkdiag}\left\{\frac{\partial \rho_{if}}{\partial \hat{x}_{11}},\dots,\frac{\partial \rho_{if}}{\partial \hat{x}_{n1}}\right\}\nm
	\\ \nm=&{\rm blkdiag}\left\{{{(\hat{x}_{11}-x_{11})^{\top}},\ldots, {(\hat{x}_{n1}-x_{n1})^{\top}}}\right\}\in\mathbb{R}^{n\times
		3n},\nm 
	\\ \nm
	{J_{f}^{\dagger}}=&\mathrm{blkdiag}\left\{\frac{(x_{11}-\hat{x}_{11})}{ \|x_{11}-\hat{x}_{11}\|^{2}} , \dots,\frac{(x_{n1}-\hat{x}_{n1})}{ \|x_{n1}-\hat{x}_{n1}\|^{2}}\right\}\in\mathbb{R}^{3n\times
		n}. \nm
	\end{align}
	Then, the desired cooperative behavior velocity of the agent $i$ is designed as 
	\begin{align}\label{v2}
	\dot{x}_{if}=&J^{\dagger}_{if} [\lambda_{f}\alpha_{if}(\tilde{\rho}_{if})-J_{i\hat{f}}\dot{\hat{x}}_{i1}],\nm
	\\
	\dot{x}_{f}=&[\dot{x}_{1f}^{\top},\dots,\dot{x}_{nf}^{\top}]^{\top}=J^{\dagger}_{f} [\Lambda_{f}\alpha_{f}(\tilde{\rho}_{f})-J_{\hat{f}}\dot{\hat{x}}_{1}],
	\end{align}
	where the gain ${\Lambda_f} : =\lambda_fI\in \mathbb{R}^{3n \times 3n}$, $\dot{\hat{x}}_{1}: =[\dot{\hat{x}}_{11}^{\top},\dots,\dot{\hat{x}}_{n1}^{\top}]^{\top} $, and $\alpha_{f}(\tilde{\rho}_{f})$ is continuous and differentiable function defined: 
	\begin{align}\label{if}
	&\alpha_{f}(\tilde{\rho}_{f}):=[\alpha_{1f}(\tilde{\rho}_{1f}),\dots,\alpha_{nf}(\tilde{\rho}_{nf})]^{\top},\\
	&\alpha_{if}(\tilde{\rho}_{if}):= \nm \\ \nm
	&\begin{cases}
	\wp_{1}\tilde{\rho}_{if}+\wp_{2}\tilde{\rho}_{if}^{[2]}, &\mathrm{if}\ \sigma_{1,i} (\dot{\tilde{\rho}}_{if},\tilde{\rho}_{if})\neq0,\ |\tilde{\rho}_{if}| \leq \phi_{s},\\
	(\beta_{1}\tilde{\rho}_{if}^{[r_{1}]}+\beta_{2}\tilde{\rho}_{if}^{[r_{2}]})^{[r_{0}]},&\mathrm{otherwise},
	\end{cases}\\
	&\sigma_{1,i_{j}} (\dot{\tilde{\rho}}_{if},\tilde{\rho}_{if}):=\dot{\tilde{\rho}}_{if}+c_{0}(\beta_{1}\tilde{\rho}_{if}^{[r_{1}]}+\beta_{2}\tilde{\rho}_{if}^{[r_{2}]})^{[r_{0}]},\nm
	\end{align}
	for $i=1,\dots,n$, where $ \wp_{1} $ and $ \wp_{2} $ are defined in \eqref{8} and \eqref{9}, and $\beta_{1},\beta_{2}, c_0, \phi_{s}, r_{0},r_{1},r_{2} \in \bR_{>0} $ are the same design parameters as in \eqref{6}.
	With the designed desired velocity 
	\eqref{v2}, the following lemma implies that each agent can maintain a desired offset with respect to the position of the
	virtual leader. 
	\begin{lemma}\label{l4}
		Consider the cooperative tracking behavior function \eqref{ctbf}. If each agent $i$ is driven by the desired velocity \eqref{v2}, then there exists a globally bounded settling time $T_{f} \in \mathbb{R}_{>0}$ such that $ \|x_{i1}-\hat{x}_{i1}\|\ge \sqrt{d_{i0}^{2}-2\phi_{s}} $, for any initial condition $\tilde{\rho}_{f}(0)$, $0 < \phi_{s} \leq \frac{d_{i0}^{2}}{2}$, and for all $i \in V$, $t\geq T_{f}$.
	\end{lemma}
	
	{We have designed a cooperative task for the flexible formation of networked agents, which guarantees the relative distance between each agent and the virtual leader, rather than the relative position between them. This means that each agent moves on a hypersphere of the leader with the radius $d_{i0}$, leading to flexibility in the formed formation. However, our design method can also be extended to solve a fixed formation problem, i.e., maintaining certain relative position among the agents. To this end, we need to steer each agent to a relative position with respect to the virtual leader. Then, we redesign the CTB function in \eqref{ctbf} and task error for each agent $i$ as
		\begin{align}\label{ctbfv}
		\rho_{if} :=\frac{1}{2} \|x_{i}-\hat{x}_{i1}-x_{ri1}\|^2, \ 
		\tilde{\rho}_{if} = - \rho_{if},
		% \label{25vv}
		%\rho_{fd}&=[\rho_{1fd}^{\top},\dots,\rho_{nfd}^{\top}]^{\top},\ \rho_{ifd}:=[0,0,0]^{\top},
		\end{align}
		where $ x_{ri1}(t) \in \bR^{3} $ is the relative position between agent $ i $ and the virtual leader. Furthermore, the corresponding Jacobian matrix $ J_{f} \in\mathbb{R}^{n\times
			3n}$ and the pseudo-inverse $J_{f}^{\dagger} \in\mathbb{R}^{3n\times n}$ become:
		{\small\begin{align*}
			J_f&={\rm blkdiag}\left\{{(x_{11}-\hat{x}_{11}-x_{ri1})^{\top}}, 
			\ldots, {(x_{n1}-\hat{x}_{n1}-x_{ri1})^{\top}}\right\}, \\
			J_{f}^{\dagger}&=\mathrm{blkdiag}
			\left\{\frac{(x_{11}-\hat{x}_{11}-x_{ri1})}{\|x_{11}-\hat{x}_{11}-x_{ri1}\|^2}, \dots,\frac{(x_{n1}-\hat{x}_{n1}-x_{ri1})}{ \| x_{n1}-\hat{x}_{n1}-x_{ri1} \|^2}\right\}.
			\end{align*}}%
		It is worth mentioning that the behavior function \eqref{ctbfv} can be also used to tackle the triangle-based 2D formation problem in \cite{ahmad2014behavioral}. For example, to achieve a triangle formation of 6 agents in Fig.~\ref{fig:01}, we can use \eqref{ctbfv} with $ x_{r11}=[0,d,0]^{\top} $, $ x_{r21}=[-\frac{\sqrt{3}d}{4},\frac{d}{4},0]^{\top} $, $x_{r31}=[\frac{\sqrt{3}d}{4},\frac{d}{4},0]^{\top} $, $ x_{r41}=[-\frac{\sqrt{3}d}{2},-\frac{d}{2},0]^{\top} $, $ x_{r51}=[0,-\frac{d}{2},0]^{\top} $, $ x_{r61}=[\frac{\sqrt{3}d}{2},-\frac{d}{2},0]^{\top} $, where $ d $ is the distance from the center of the equilateral triangle to its three vertices, i.e., 1, 4, and 6.
		\begin{figure}[!tp]\centering
			\includegraphics[width=3.5cm,trim=5 5 5 2, clip]{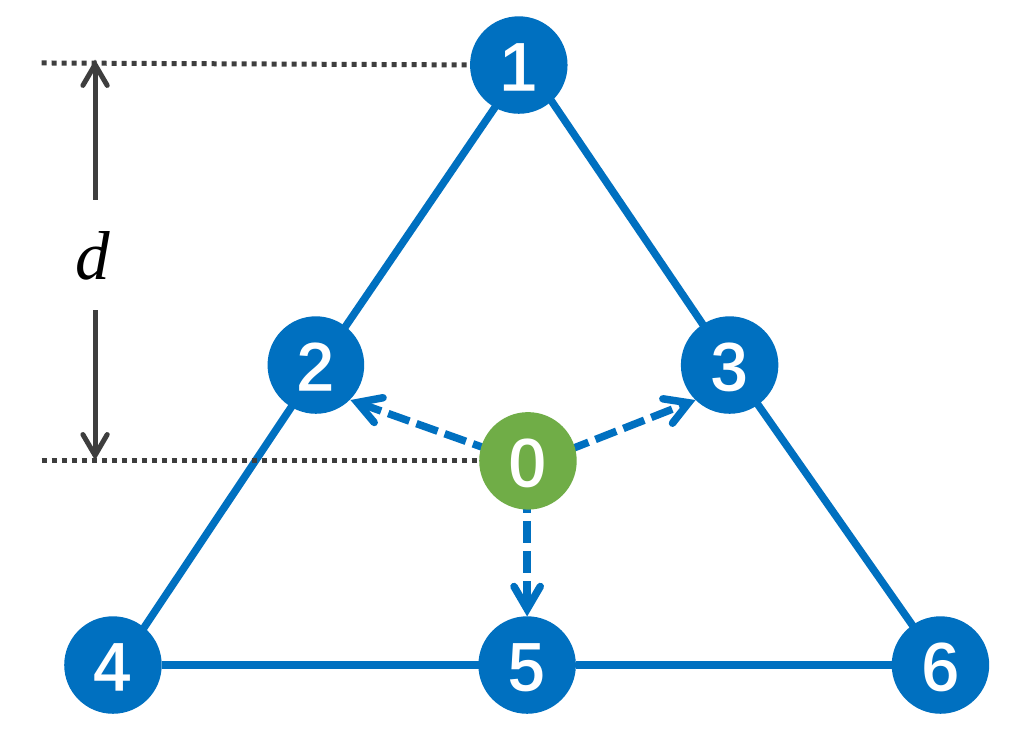}
			\caption{Triangle formation.}
			\label{fig:01}
		\end{figure}
	}
	{\begin{remark}
			Note that the estimator \eqref{11} can not be directly extended to the directed communication topology due to the assumption of a positive definite $ H $. If the estimator \eqref{11} can be redesigned to be applicable to the digraph case, then the rest of the results in this paper are still valid.
	\end{remark}}\subsection{Merged Desired Velocities}
\begin{figure}
		\centering\includegraphics[width=6.5cm,trim=0 0 0 0, clip]{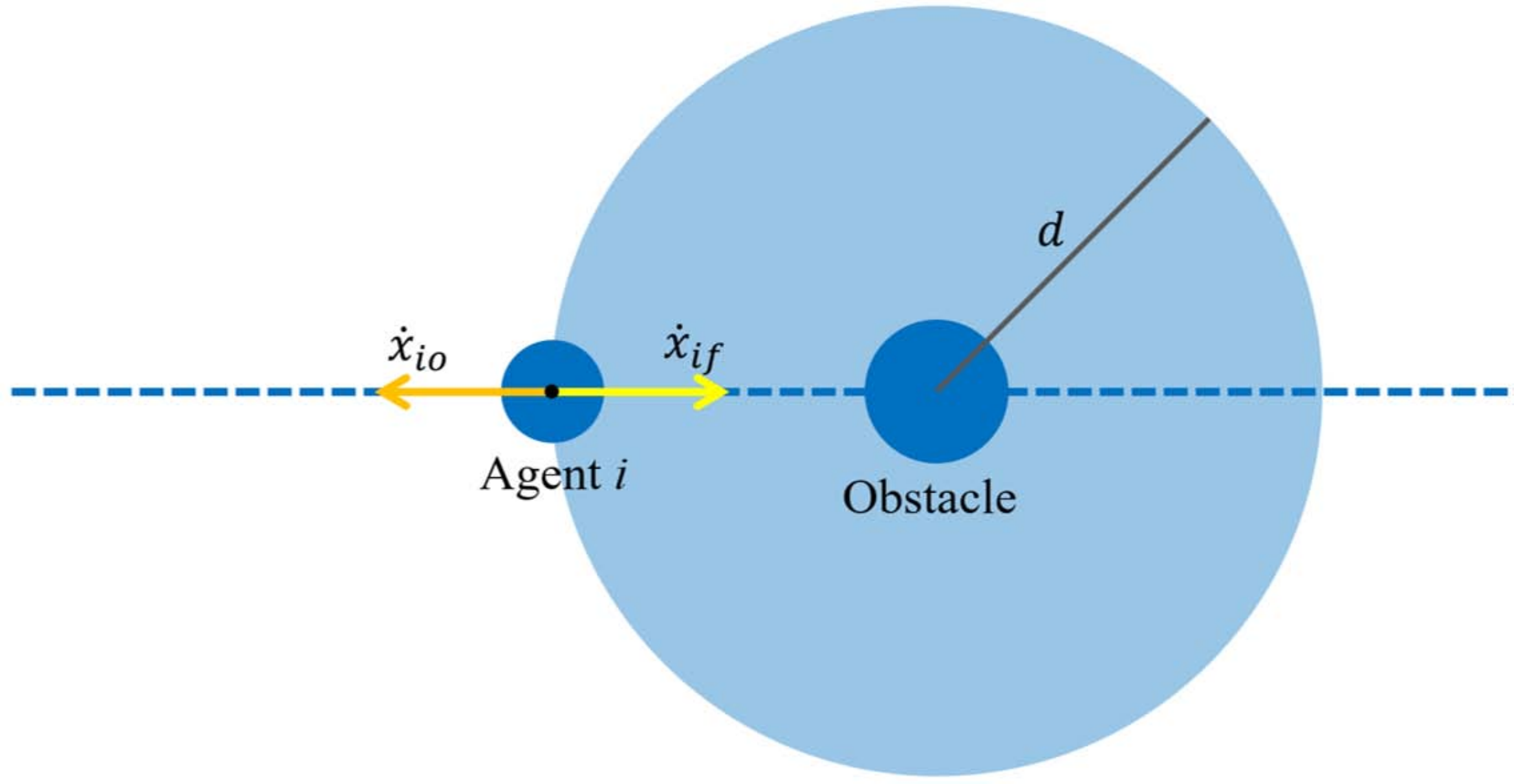}
		\caption{The illustration of the {local minima}, where  $\dot{x}_{io} = \dot{x}_{if}$.}
		\label{fig.0}
	\end{figure}

	In this section, we design the {desired velocity guaranteeing fixed-time convergence} for each agent when the two behaviors, namely, the collision avoidance behavior and the cooperative tracking behavior are combined. This combination is taken using the null-space-based behavioral approach, which essentially prioritizes the two tasks. 
	%	The higher-priority task is consistently carried out while the lower-priority task is projected onto the null space of the higher-priority task such that it is executed when it does not oppose the execution of the higher-priority task. 
	{The challenge of the combination is to solve the {local minima} when {$ \dot{x}_{io}=-\dot{x}_{if} $}. The {local minima} is illustrated  in Fig.~\ref{fig.0}, where $\dot{x}_{io}$ and $\dot{x}_{if}$ are colinear with the same magnitudes but opposite directions. We solve this problem by designing a robust term in the combination at the {local minima}.}
	
	To guarantee the security of system operation, the obstacle/collision avoidance behavior is given a higher priority in this paper. Then,
	with the framework of NSB, we determine the desired velocity $\dot{x}_{d}:\bR_{>0}\to\bR^{3n}$ of the agent $ i $ by merging the two designed velocities \eqref{v1} and \eqref{v2}, which leads to
	{\begin{align}\label{v}
		\dot{x}_{id}=&\begin{cases}
		\dot{x}_{io}+(I-J_{io}^{\dag}J_{io})\dot{x}_{if}, &\mathrm{if}\ {\|\dot{x}_{io}+\dot{x}_{if}\|\neq 0 ,}\\
		\dot{x}_{io}+(I-J_{io}^{\dag}J_{io})\dot{x}_{if}+\dot{x}^{c}_{id},&\mathrm{if}\ {\|\dot{x}_{io}+\dot{x}_{if}\|= 0 ,}
		\end{cases}
		\end{align}
		$ \forall~i \in V $. The additional term
		\begin{align}\label{vc}
		\dot{x}^{c}_{id}=\int_{t}^{t+1}\delta_{d}{\cal{R}}_{z}(\theta_{z}){\cal{R}}_{y}(\theta_{y}){\cal{R}}_{x}(\theta_{x})\frac{x_{i1}-x_{i}^{o}}{\|x_{i1}-x_{i}^{o}\|}dt
		\end{align}
		is added only when the {local minima} $\|\dot{x}_{io}+\dot{x}_{if}\|=0 $ is reached,  where
		$ \delta_{d} \in \mathbb{R}_{>0} $ is a small constant, $ \theta_{x},\theta_{y},\theta_{z}\in[-\frac{\pi}{2},\frac{\pi}{2}] $ are random or pre-determined by the users that satisfy $ \theta_{x}^{2}+\theta_{y}^{2}+\theta_{z}^{2}\neq 0 $, and $ {\cal{R}}_{x}(\theta_{x}), {\cal{R}}_{y}(\theta_{y}) , {\cal{R}}_{z}(\theta_{z}) \in\bR^{3\times 3}$ are three rotation matrices that rotate the  vector $ \mathrm{sgn}{(x_{i1}-x_{i}^{o})}$  by $ \theta_{x},\theta_{y}$ and $\theta_{z} $ around $ x$-, $ y$- and $ z$-axis, respectively. {Without the robust term $\dot{x}^{c}_{id}$, the desired velocity $\dot{x}_{id}$ may cause the agent to stuck at the {local minima}.}}
	
	The velocities $ \dot{x}_{id} $, can be formed as a stacked column vector: $\dot{x}_{d}=[\dot{x}_{1d}^{\top},\dots,\dot{x}_{nd}^{\top}]^{\top}$, from which we further compute the position references ${x}_{d}$ by applying the closed-loop inverse kinematics algorithm \cite{antonelli2009stability}. With the merged velocity 
	\eqref{v}, we prove that each agent in the network can achieve the both tasks simultaneously with a fixed settling time.
	\begin{theorem}\label{taskl4}
		Consider the obstacle avoidance behavior in \eqref{1} and the cooperative behavior in \eqref{ctbf}. If each agent is driven by the merged desired velocity in \eqref{v}, then there exists a settling time $T_{i}$ such that 
		$$
		\|x_{i1}-x_{i}^{o}\|\ge {\tilde{d}}, \ \text{and} \ \|x_{i1}-\hat{x}_{i1}\|\ge \sqrt{d_{i0}^{2}-2\phi_{s}}
		$$
		for any $(\tilde{\rho}_{io}(0),\tilde{\rho}_{if}(0))\in\bR\times\bR$, $0 < \phi_{s} \leq \min \{d^2/2, d_{i0}^{2}/2\}$, and for all $i \in V$ and {$t\geq T_{i}$, where $ \tilde{d}:=\sqrt{d^{2}-2\phi_{s}} $.} 
	\end{theorem}
	Note that the proof of the theorem is not just a simple combination of the conclusions of Lemma~\ref{l1} and Lemma~\ref{l4}. When there is no conflict between the two tasks, we have $J_f J_{io}^{\dag} = 0$, which means that two tasks in the velocity space are orthogonal and thus the fixed-time properties can be proved independently. However, if the tasks are conflicting, i.e., $J_f J_{io}^{\dag} \ne 0$, then the proof becomes nontrivial, {which is more challenging than the case $J_f J_{io}^{\dag} = 0$}. We present the detailed proof in Appendix~\ref{Ap:task14}.
	\begin{remark}
		The designed velocity in \eqref{v} is inspired by the works in \cite{chiaverini1997singularity,antonelli2006kinematic}, which also provide another solution for the {local minima} problem by intentionally adding measurement noises. However, this method does not provide a theoretical guarantee for the convergence of the behavior errors in \eqref{2} and \eqref{15}.
		Different from \cite{chiaverini1997singularity,antonelli2006kinematic}, the robust term $ \dot{x}^{c}_{id} $  in \eqref{v} is integrable and differentiable, which is essential for   guaranteeing the fixed-time convergence of the tracking errors. 
	\end{remark}{\begin{remark}
			The proposed approach \eqref{v} can be extended to the multiple-task case following the merging rule developed in \cite{antonelli2008stability}, but the merged desired velocity at the {local minima} should be redesigned, and it is difficult to prove the fixed-time convergence of all the task errors, especially when the tasks are conflicting with each other.
	\end{remark}}
	
	\section{Robust Tracking Control Design}\label{sec4}
	{In this section, we design a robust control strategy to track the desired velocity in \eqref{v}.}
	Utilizing the final desired velocity $\dot{x}_{d}$ and position ${x}_{d}$, the two tracking errors are defined as
	\begin{align*}
	\tilde{x}_{i1}:=x_{i1}-x_{id}, \ {\tilde{x}}_{i2}:={x}_{i2}-\dot{x}_{id}.
	\end{align*}
	
	First, we design a {FTTSM}, following the design of \cite{jiang2016fixed,Zhou2019FTC}, which is constructed by two functions {$\sigma_{1i}: \mathbb R^3 \times \mathbb R^3 \to \mathbb R^3$ and $\sigma_{2i}: \mathbb R^3 \times \mathbb R^3 \to \mathbb R^3$ for each $i \in V$:
		\begin{align*}
		\sigma_{1i}(\tilde{x}_{i1},\tilde{x}_{i2})&:=\tilde{x}_{i2}+c_{1}\tilde{x}_{i1}+c_{2}(\beta_{1}\tilde{x}_{i1}^{[r_{1}]}+\beta_{2}\tilde{x}_{i1}^{[r_{2}]})^{[r_{0}]},\\
		\sigma_{2i}(\tilde{x}_{i1},\tilde{x}_{i2})&:=\tilde{x}_{i2}+c_{1}\tilde{x}_{i1}+c_{2}(\wp_{1}\tilde{x}_{i1}+\wp_{2}\tilde{x}_{i1}^{[2]}), 
		\end{align*}}%
	where $ \wp_{1} $ and $ \wp_{2} $ are designed as \eqref{8} and \eqref{9}, 
	and $c_{1}$, $c_{2}$, $\phi_{s}$, $ r_{0},r_{1},r_{2}\in \bR_{>0}$ with $r_{1}r_{0}>1$, $r_{2}r_{0}<1$, and $r_{0}\in(\frac{1}{2},1) $. 
	
	Next, by using $\sigma_{1i}$ and $\sigma_{2i}$, we define the FTTSM ${\mathcal{S}}_{i}: \mathbb R^3 \times \mathbb R^3 \to \mathbb R^3$ as
	\begin{small}
		\begin{align}\label{S}
		&{\mathcal{S}}_{i}(\cdot):=[{\mathcal{S}}_{i_{1}}(\cdot),{\mathcal{S}}_{i_{2}}(\cdot),{\mathcal{S}}_{i_{3}}(\cdot)]^{\top},\ i \in V,\ j=1,2,3,\nm\\
		&{\mathcal{S}}_{i_{j}}(\sigma_{1i_{j}},\sigma_{2i_{j}})\\
		&:=\begin{cases}
		\varrho\sigma_{2i_{j}} (\tilde{x}_{i1_{j}},\tilde{x}_{i2_{j}}), &\mathrm{if}\ \sigma_{1i_{j}} (\tilde{x}_{i1_{j}},\tilde{x}_{i2_{j}})\neq0,\ |\tilde{x}_{i1_{j}}| \leq \phi_{s},\\
		\varrho\sigma_{1i_{j}} (\tilde{x}_{i1_{j}},\tilde{x}_{i2_{j}}),&\mathrm{otherwise},
		\end{cases} \nm
		\end{align}
	\end{small}%
	{which is a piecewise continuous and differentiable function with the design parameter $ \varrho\in \bR_{>0} $. }
	For the FTTSM ${\mathcal{S}}_{i}$, we have the following result {in Lemma \ref{l7}}. 
	\begin{lemma}\label{l7}
		Consider the FTTSM $ {\mathcal{S}}_{i}$ defined by \eqref{S}. For any $ \delta_{s}>0 $, $ \tilde{x}_{i1_{j}}(0),\ \tilde{x}_{i2_{j}}(0)\in\bR $, if $ |{\mathcal{S}}_{i_{j}}(\sigma_{1i_{j}},\sigma_{2i_{j}})|\le\delta_{s} $, then there exists a fixed time $ T_{s_{*}} >0$ such that
		\begin{align*}
		|\tilde{x}_{i1_{j}}(t)|\le\max \{\phi_{s},\delta_{s_{1}}\},\ |\tilde{x}_{i2_{j}}(t)|\le \max\{\delta_{s_{0}},\delta_{s_{2}}\},
		\end{align*} 
		for all $ t\ge T_{s_{*}} $, where
		\begin{align*}
		\delta_{s_{0}}& := \delta_{s}/\varrho+ c_{1}\phi_{s}+c_{2}(\beta_{1}\phi_{s}^{r_{1}}+\beta_{2}\phi_{s}^{r_{2}})^{r_{0}} ,\\
		\delta_{s_{1}} &:= \min\left\{\frac{\delta_{s}}{\varrho \tilde{c}_{1}},
		\frac{(\delta_{s}/\varrho \tilde{c}_{2})^{1/r_{0}r_{1}}}{(\tilde{\beta}_{1})^{1/r_{1}}},
		\frac{(\delta_{s}/ \varrho \tilde{c}_{2})^{1/r_{0}r_{2}}}{(\tilde{\beta}_{2})^{1/r_{2}}}\right\} , 
		\\
		\delta_{s_{2}} &:= \delta_{s}/ \varrho+ c_{1}\delta_{s_{1}}+c_{2}(\beta_{1}\delta_{s_{1}}^{[r_{1}]}+\beta_{2}\delta_{s_{1}}^{[r_{2}]})^{[r_{0}]}.
		\end{align*}
		In the definition of $\delta_{s_{1}}$, the parameters $\tilde{c}_{1}$,  $\tilde{c}_{2}$, $\tilde{\beta}_{1}$, $\tilde{\beta}_{2} \in \mathbb{R}_{>0}$ are selected such that
		$\tilde{c}_{1} \leq c_1$,  $\tilde{c}_{2} \leq c_2$, $\tilde{\beta}_{1} \leq \beta_1$, $\tilde{\beta}_{2} \leq \beta_2$, and at least one of the following inequalities holds.
		\begin{subequations} \label{eq:tildecbeta}
			\begin{align}
			|\tilde{x}_{i1_{j}}(0)| \ge {\delta_{s}}/{(\varrho \tilde{c}_{1})}>0 ,
			\\
			|\tilde{x}_{i1_{j}}(0)|\ge \frac{(\delta_{s}/\varrho \tilde{c}_{2})^{1/r_{0}r_{1}}}{(\tilde{\beta}_{1})^{1/r_{1}}}>0 ,\\
			\ |\tilde{x}_{i1_{j}}(0)|\ge \frac{(\delta_{s}/ \varrho \tilde{c}_{2})^{1/r_{0}r_{2}}}{(\tilde{\beta}_{2})^{1/r_{2}}}>0.
			\end{align}
		\end{subequations}
	\end{lemma}
	{Lemma~\ref{l7} essentially implies that a norm-bounded ${\mathcal{S}_{i}}$ guarantees the boundedness of $\tilde{x}_{i1}$ and $\tilde{x}_{i2}$ in their magnitudes. Then, our design problem becomes how to design a controller to stabilize ${\mathcal{S}_{i}}$. 
		
		First, to guarantee the fixed-time convergence of velocity and position tracking errors in the reaching phase, we present a reaching law as 
		\begin{align}
		u_{i}^{S}(t)=&-k_{1i}(t){\mathcal{S}_{i}}^{[\gamma_{1}]}-{k}_{2i}(t){\mathcal{S}_{i}}^{[\gamma_{2}]},
		\label{eq:uS}\\
		\label{k1} k_{1i}(t)=&k_{1i}^{0}+\frac{k_{1i}^{M}-k_{1i}^{0}}{\exp[-c_{s1}(t-T_{0})]+1}, \\ {k}_{2i}(t)=&k_{2i}^{0}+\frac{k_{2i}^{M}-k_{2i}^{0}}{\exp[-c_{s2}(t-T_{0})]+1}, \label{k2}
		\end{align}where $ \gamma_{1}>1$, $0<\gamma_{2}<1 $, $ c_{s1}, c_{s2}, T_{0}\in\bR_{>0} $. The two gain functions $ k_{1i}$ and $k_{2i}$ are bounded in the intervals $[k_{1i}^{m},k_{1i}^{M})$ and $[k_{2i}^{m},k_{2i}^{M})$, respectively, where $ k_{1i}^{M}>k_{1i}^{0}>0 $, $ k_{2i}^{M}>k_{2i}^{0}>0 $, and 
		\begin{align*}
		k_{1i}^{m}= k_{1i}^{0}+\frac{k_{1i}^{M}-k_{1i}^{0}}{\exp(c_{s1}T_{0})+1}, \ k_{2i}^{m}= k_{2i}^{0}+\frac{k_{2i}^{M}-k_{2i}^{0}}{\exp(c_{s2}T_{0})+1}.
		\end{align*}  
		The functionality of the parameters in $ k_{1i} (t)$ and $ k_{2i} (t)$ is discussed. 
		The scalars $c_{s1}$ ($c_{s2}$) determines the slope of $ k_{1i} (t)$ ($ k_{2i} (t)$) at time $t=T_0$. A larger value of $c_{s1}$ ($c_{s2}$) implies a steeper slope. The value of $k_{1i}^{M}$ ($k_{1i}^{M}$) determines the upper bound of the gain $ k_{1i}$ (($k_{2i}$)).
		Observe that the gains $ k_{1i} (t)$ and $ k_{2i} (t)$  are designed as monotonically increasing and bounded functions. The purpose of using these time-varying gains is to improve the tracking accuracy  in the reaching phase. The parameter $ T_{0} $ in \eqref{k1} and \eqref{k2} is chosen to be larger than the convergence time of the tracking errors, namely 
		\begin{align*}
		T_{0}>T_{0}^{*} :=  T_{i}(\tilde{\rho}_{if}(0))+T_{S1}+T_{s}.
		\end{align*}
		In Remark~\ref{rem5}, we provide further discussion on designing $k_{1i}(t)$ and ${k}_{2i}(t)$.
		
		Second, to address the uncertainties in the system model \eqref{agent} and the derivative of sliding-mode $\mathcal{S}_i$ in \eqref{S}, we define the  function $ F_{i}(z_{i}) $ for the agent $i$ to collect these uncertainties:
		\begin{align}\label{32}
		F_{i}(z_{i}):=&f_{i}(\bar{x}_{i})+\chi_{i}(z_{1i}),\\ \nm
		z_{i}:=&[{x}_{i1}^{\top},{x}_{i2}^{\top},\ddot{x}_{id}^{\top},\dot{\tilde{x}}_{i1}^{\top},\dot{\alpha}_{i}(\dot{\tilde{x}}_{i1})^{\top}]^{\top},\\ \chi_{i}(z_{1i}):=&-\ddot{x}_{id}+c_{1}\dot{\tilde{x}}_{i1}+c_{2}\dot{\alpha}_{i}(\dot{\tilde{x}}_{i1}),\nm
		\end{align}
		where ${\alpha}_{i}(\tilde{x}_{i1}):=[{\alpha}_{i_{1}}(\tilde{x}_{i1_{1}}),{\alpha}_{i_{2}}(\tilde{x}_{i1_{2}}),{\alpha}_{i_{3}}(\tilde{x}_{i1_{3}})]^{\top}$, with 
		\begin{align*}
		&{\alpha}_{i_{j}}(\tilde{x}_{i1_{j}}):=\\
		&\begin{cases}
		(\beta_{1}\tilde{x}_{i1_{j}}^{[r_{1}]}+\beta_{2}\tilde{x}_{i1_{j}}^{[r_{2}]})^{[r_{0}]}, &\mathrm{if}\ \sigma_{1i_{j}} (\tilde{x}_{i1_{j}},\tilde{x}_{i2_{j}})\neq0,\ |\tilde{x}_{i1_{j}}| \leq \phi_{s},\\\wp_{1}\tilde{x}_{i1_{j}}+\wp_{2}\tilde{x}_{i1_{j}}^{[2]},&\mathrm{otherwise}.
		\end{cases}\nm
		\end{align*}
		for any $ i \in V$ and $ j=1,2,3 $. The parameters $\beta_{1}$, $\beta_{2}$, $\wp_{1}$, and $\wp_{2}$ are also used in \eqref{6} and 
		\eqref{if}. Then, we apply the universal approximation property of RBFNNs to estimate the uncertainty function $F_{i}(z_{i})$ by {an} artificial neuronal network. We refer to e.g., \cite{chen2014adaptive,Zhou2019FTC} for more details. Specifically, for any $ \varepsilon_{Ni} \in \mathbb{R}_{>0}$, we find a RBFNN with $ h_{i} $ neurons such that $F_{i}(z_{i})$ is described by
		\begin{align}\label{33}
		F_{i}(z_{i}):=&{W_{i}^{\ast}(t)}^{\top}\phi_{i}(z_{i})+\varepsilon_{i},\ \forall t\in\bR_{\ge 0},z_{i}\in\Omega_{zi},
		\end{align}
		with $ \|\varepsilon_{i}\|\le\varepsilon_{Ni} $, $\Omega_{zi}\subset\bR^{15}$ a prefixed sufficiently large compact set, and $ {W_{i}^{\ast}(t)}\in\bR^{h_{i}\times 3} $ the weight matrix. The basis function $ \phi_{i}:\bR^{15} \rightarrow \bR^{h_{i}} $ is defined with its $ k $th element as
		$$\phi_{ik}(z_{i}):=\exp[-(z_{i}-\mu_{ik})^{\top}(z_{i}-\mu_{ik})/\psi_{ik}^2],\ {k=1,\dots,h_{i}},$$
		where $\mu_{ik}\in \mathbb{R}^{15}$ denotes the center of the receptive field, and $\psi_{ik}$ is the width of the Gaussian function.  Using \eqref{33}, we then construct an intelligent controller $ u_{i}^{N} (t)$ for the uncertainty compensation:
		\begin{align}\label{ad1}
		u_{i}^N(t)= -\hat{W}_{i}^{\top}\phi_{i}(z_{i}),\ \text{with} \
		\dot{\hat{W}}_{i}(t)= \Gamma_{i}\phi_{i}(z_{i}){\mathcal{S}_{i}}^{\top},
		\end{align}where $\Gamma_{i}{\in\bR^{h_{i}\times h_{i}}}$ is a positive definite
		constant gain matrix, and $\hat{W}_{i}$ denotes the estimation of $W_{i}^{\ast }$.

		Third, to reject the external disturbance $ d_{i} $ and the internal disturbance $ \varepsilon_{i} $ generated from \eqref{33}, we design a disturbance compensator $ u_{i}^{C}(t) $ as 
		\begin{align} \label{ad2}
		u_{i}^C(t)= -\hat{\delta}_{i}\text{sign}({\mathcal{S}_{i}}), \ \text{with} \
		\dot{\hat{\delta}}_{i}(t)= \gamma_{3i}\|{\mathcal{S}_{i}}\|_{1},
		\end{align}
		where $\gamma_{3i} \in \mathbb{R}_{>0}$ is a design parameter, and $\hat{\delta}_{i}$ represents the upper bound of the disturbances $\delta_{i}$, namely, $\|\varepsilon_{i}+d_{i}\|\leq \delta_{i}$.

		Finally, we combine the reaching law $ u_{i}^{S} (t) $, the intelligent controller $ u_{i}^{N} (t) $ and the disturbance compensator $ u_{i}^C(t) $ to provide a fixed-time behavioral control law  as follows:
		\begin{align}\label{u}
		u_{i}(t)=u_{i}^{S}(t)+u_{i}^N(t)+u_{i}^C(t).
		\end{align}

		Before proceeding to {theoretical} analysis of the above control law, we present {a lower bound} of the design parameter $\lambda_{io}$ in \eqref{v1}.  
		Note that {the sliding mode \eqref{S} contains both position and velocity errors}, and the NSB approach is a kinematic acting on the second-order dynamics \eqref{agent} through the desired velocity \eqref{v} rather than a desired position. It is required that the velocity error dominates the position error in the sliding manifold \eqref{S}. More precisely, $ \|x_{i1}-x_{i}^{o}\|\le d $ should be guaranteed, i.e., the following inequality should hold.}
	\begin{align}\label{25}
	&(\dot{x}_{i1}-J_{io}^{\dag}\lambda_{io}\tilde{\rho}_{io})^{\top}(\dot{x}_{i1}
	-J_{io}^{\dag}\lambda_{io}\tilde{\rho}_{io})>\\
	&c_{1}^{2}\tilde{x}_{i1}^{\top}\tilde{x}_{i1}+c_{2}^{2}\alpha_{i}(\tilde{x}_{i1})^{\top}
	\alpha_{i}(\tilde{x}_{i1})+2c_{1}(\dot{x}_{i1}-J_{io}^{\dag}\lambda_{io}\tilde{\rho}_{io})^{\top}\tilde{x}_{i1} \nm\\&+2c_{1}
	c_{2}\tilde{x}_{i1}^{\top}
	\alpha_{i}(\tilde{x}_{i1})
	+2c_{2}(\dot{x}_{i1}
	-J_{io}^{\dag}\lambda_{io}\tilde{\rho}_{io})^{\top}\alpha_{i}(\tilde{x}_{i1}). \nm
	\end{align}
	{Taking the norm on both sides of \eqref{25} then leads to}
	\begin{align}\label{Constraint}
	\lambda_{io}> \frac{-\Upsilon_{2i}+\sqrt{\Upsilon_{2i}^{2}-4\Upsilon_{1i}\Upsilon_{3i}}}{2\Upsilon_{1i}}: = \lambda^{\ast}_{i},
	\end{align}
	where  $\Upsilon_{1i}=\|J_{io}^{\dag}\|^{2}\alpha_{io}^{2}$,
	$\Upsilon_{2i}=-2\|J_{io}^{\dag}\||\alpha_{io}|(\|\dot{x}_{i1}\|+c_{1}\|\tilde{x}_{i1}\|+c_{2}\|\alpha_{si}(\tilde{x}_{i1})\|)$, and
	$\Upsilon_{3i}=-(\|\dot{x}_{i1}\|^{2}+c_{1}^{2}\|\tilde{x}_{i1}\|^{2}+c_{2}^{2}\|\alpha_{si}(\tilde{x}_{i1})\|^{2}
	+2c_{1}\|\dot{x}_{i1}\|\|\tilde{x}_{i1}\|+2c_{2}\|\dot{x}_{i1}\|\|\alpha_{si}(\tilde{x}_{i1})\|+2c_{1}c_{2}
	\|\tilde{x}_{i1}\|\|\alpha_{si}(\tilde{x}_{i1})\|)$. {The inequality} \eqref{Constraint} can be recognized as a constraint on the state-dependent gain $ \lambda_{io} $, and we choose $\lambda_{io}=
	\lambda^{\ast}_{i}+\varrho_{i}$, {with $\varrho_{i} \in \mathbb{R}_{>0}$ an arbitrary small robust term, to make the closed-loop system robust to measurement noises.
		
		Now, we are ready to provide a theoretical guarantee on the performance of the control law \eqref{u}. The detailed proof is given in Appendix~\ref{Ap:the}.}
	\begin{theorem}\label{the1}
		Consider $n$ networked autonomous agents with second-order
		nonlinear dynamics in~\eqref{agent}, and each of the agents are controlled by \eqref{u}. Suppose that \eqref{a3.1} holds, and the design parameter $ \lambda_{io} $ is selected to {fulfill} $ \lambda_{io}>\max\{\lambda_{io1},\lambda_{io2},\lambda^{\ast}_{i}\} $. Then, for a given ${\delta}_{s} \in \mathbb{R}_{>0}$, there exists a finite  $T_{S1} \in \mathbb{R}_{>0}$ such that $ \|{\mathcal{S}}_{i}\|\le \delta_{s} $ holds for all $t \ge T_{S1}$.
	\end{theorem}
	{\begin{remark}
			\label{rem:parameter}
			The selection of the parameters in the controller \eqref{u} is discussed.  Generally, when tuning these parameters, there is a trade-off between the convergence speed and track precision. Specifically, the larger values of $ k_{1i}^{0} $, $ k_{2i}^{0} $, $ k_{1i}^{M} $, $ k_{2i}^{M} $, and $ \gamma_{1} $, or the smaller values of $ \gamma_{2} $, will lead to a faster convergence speed. Moreover, the smaller values of $\phi _s$, $c_1$, $c_2$, $\beta _1$, $\beta _2$, and $r_2r_{0}$, or larger values of $ \varrho $, $r_1r_{0}$, yield a higher convergence precision. 
	\end{remark}}
	\section{Design of Adjustable Control Gains\label{rem5}}
	
	An alternative design of the gains $ k_{1i} (t)$ and $ k_{2i} (t)$ in \eqref{eq:uS} can be considered. Note that the gains in \eqref{k1} and \eqref{k2} can only be increased. Such a design may not be enough to meet the requirements of high-precision formation control in e.g., 3D synthesis imaging using networked vehicles \cite{Negahdaripour2009Opti}.
	In this case, we can use $ k_{1i} (t)$ and $ k_{2i} (t)$ in a more general form as
	\begin{align}\label{k1v}
	k_{1i}(t):=&k_{1i}^{0}+\sum_{j=1}^{\iota}\frac{(-1)^{\zeta_{ji0}}\zeta_{ji}(k_{j1i}^{M}-k_{1i}^{0})}{\exp[-c_{ji1}(t-T_{j0})]+1},\\ {k}_{2i}(t):=&k_{2i}^{0}+\sum_{j=1}^{\iota}\frac{(-1)^{\zeta_{ji0}}\zeta_{ji}(k_{j2i}^{M}-k_{2i}^{0})}{\exp[-c_{ji2}(t-T_{j0})]+1},\label{k2v}
	\end{align}
	where $ \iota\in\bN_{\ge1} $, $ c_{ji1}, c_{ji2}  \in \mathbb{R}_{>0}$, $ T_{j0}\in\bR_{>0} $ and $ T_{10}<T_{20}< \dots<T_{\iota 0}$. The two gains are bounded in the intervals $[\kappa_{1i}^{m},\max_{j}\{k_{ j1i}^{M}\})$ and $[\kappa_{2i}^{m},\max_{j}\{k_{j 2i}^{M}\})$, respectively, where
	\begin{align*}
	\kappa_{1i}^{m}&= k_{1i}^{0}+\sum_{j=1}^{\iota}\frac{(-1)^{\zeta_{ji0}}\zeta_{ji}(k_{j1i}^{M}-k_{1i}^{0})}{\exp(c_{js1}T_{j0})+1}\\ \kappa_{2i}^{m}&= k_{2i}^{0}+\sum_{j=1}^{\iota}\frac{(-1)^{\zeta_{ji0}}\zeta_{ji}(k_{j2i}^{M}-k_{2i}^{0})}{\exp(c_{js2}T_{j0})+1}.
	\end{align*}
	The binary parameters $\zeta_{ji}$ and $\zeta_{ji0}$ are valued such that $\zeta_{ji} = 0$ yields constant {gains} $k_{1i}(t) = k_{1i}^{0}$ and $k_{2i}(t) = k_{2i}^{0}$; $\zeta_{ji} = 1$ provides a change of gains; $\zeta_{ji0} = 1$ leads to the increase of the gains, while $\zeta_{ji0} = 0$ means to decrease the gain. To ensure positive gains $ k_{1i}(t), k_{2i}(t)$, i.e., $\kappa_{1i}^{m}, \kappa_{2i}^{m} \in \mathbb{R}_{>0}$, we impose  the constraints:  (1) If $ \zeta_{ji0}=1 $, $ k_{j1i}^{M}<2 k_{1i}^{0} $ and $ k_{j2i}^{M}<2 k_{2i}^{0} $, for all $j = 1, 2, \cdots, \iota$. (2) If $ \zeta_{ji0}=0 $, $ k_{\iota1i}^{M}>\dots>k_{21i}^{M}>k_{11i}^{M}>k_{1i}^{0}$ and $ k_{\iota2i}^{M}>\dots>k_{22i}^{M}>k_{12i}^{M}>k_{2i}^{0}$.
	
\begin{figure}[!tp]\centering
		\includegraphics[width=7.5cm,trim=110 280 130 295, clip]{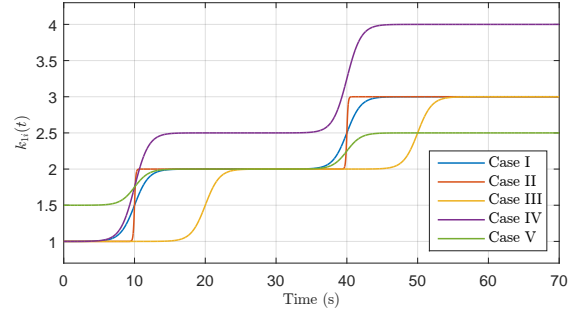}
		\caption{Function $ k_{1i} (t)$ in (\ref{k1v}) with $ \zeta_{ji0}=0 $ and $ \zeta_{ji}=1 $.}
		\label{fig.1}
	\end{figure}

\begin{figure}[!tp]\centering
		\includegraphics[width=7.5cm,trim=110 280 130 295, clip]{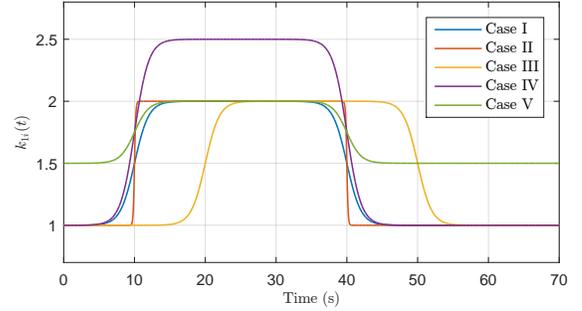}
		\caption{Function $ k_{1i} (t)$ in (\ref{k1v}) with $ \zeta_{1i0}=0 $ and $ \zeta_{2i0}=\zeta_{ji}=1 $.}
		\label{fig.1v}
	\end{figure}

	The generalizations in \eqref{k1v} and \eqref{k2v} provide extra freedom to adjust the control gains in \eqref{eq:uS}, which can further lead to a higher control accuracy in the reaching phase. Moreover, with the gains in \eqref{k1v} and \eqref{k2v}, the conclusion of Theorem~\ref{the1} still holds. The proof follows similarly as the case that $ k_{1i} (t)$ and $ k_{2i} (t)$ are designed as \eqref{k1} and \eqref{k2}. Thus, we omit the details. To illustrate the gains $ k_{1i} (t)$ and $ k_{2i} (t)$, we select the design parameters of $ k_{1i} (t)$ as in Table~\ref{tab1} and visualize the function $ k_{1i} (t)$ in Figs.~\ref{fig.1} and \ref{fig.1v}.
	\begin{table}[h]
		\caption{Parameters of $ k_{1i} (t)$ in (\ref{k1v}) with $ \iota=2 $}
		\label{tab1}
		\begin{center}
			\begin{tabular}{c||c|c|c|c|c|c|c}\hline
				Cases  &$k_{1i}^{0}$& $k_{11i}^{M}$& $k_{21i}^{M}$& $c_{js1}$&$T_{10}$&$T_{20}$ &colour\\ \hline \hline
				Case I & 1& 2&3& 1& 10&40 &blue\\ \hline
				Case II & 1& 2&3& 11& 10&40& orange\\ \hline
				Case III& 1& 2&3& 1& 20&50&yellow\\ \hline
				Case IV & 1&2.5&4& 1& 10&40& purple\\ \hline
				Case V& 1.5& 2&2.5& 1& 10&40& blue\\ \hline
			\end{tabular}
		\end{center}
	\end{table}
	
	\begin{figure}[!tp]\centering
		\includegraphics[width=4cm,trim=0 0 0 0, clip]{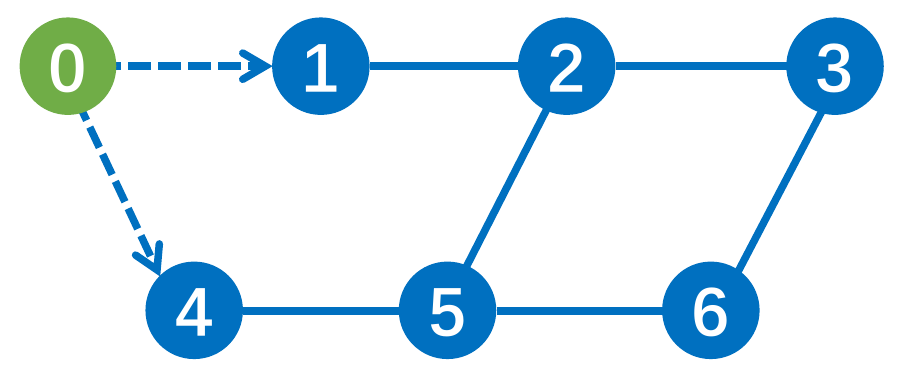}
		\caption{The autonomous agents network.} 
		\label{fig:1}
	\end{figure}

	\section{Simulation Results}\label{sec5}
	Consider a multi-agent {system} connected by an undirected network as shown in Fig.~\ref{fig:1}, which contains 6 followers and a virtual leader in a 3-dimensional space. 
	The graph $G$ associated with the communication network is unweighted, i.e., $ a_{ij}=a_{ji}=1 $ if there exists an information exchange between the agents $ i $ and $ j $, and $ b_{i}=1 $ if the agent $ i $ can obtain the information from the leader. 
	
	Each follower agent is modeled by a second-order nonlinear dynamic system as in~\eqref{agent} with nonlinear terms   \\
	$ f_{1}(\bar{{x}}_{1})=0.1\|{x}_{11}\|\sin({x}_{12}), \ f_{2}(\bar{{x}}_{2})=0.1\|{x}_{21}\|\tanh({x}_{22}),  $\\
	$ f_{3}(\bar{{x}}_{3})=0.1\|{x}_{32}\|\sin({x}_{31}),\ f_{4}(\bar{{x}}_{4})=0.1\|{x}_{42}\|\tanh({x}_{41}), $ \\
	$ f_{5}(\bar{{x}}_{5})=0.1\|{x}_{51}\|\sin({x}_{52}),\ f_{6}(\bar{{x}}_{6})=0.1\|{x}_{61}\|\sin({x}_{62}). $\\	
	and disturbances\\
	$ d_{1}(\bar{{x}}_{1},t)=0.5\|0.02{x}_{11}\|[\sin(0.5t),\sin(0.7t),\cos(0.5t)]^{\top}, $ \\  $ d_{2}(\bar{{x}}_{2},t)=0.5\|0.03{x}_{21}\|[\sin(0.5t),\tanh(0.7t),\cos(0.5t)]^{\top}, $ \\
	$ d_{3}(\bar{{x}}_{3},t)=0.5\|0.02{x}_{31}\|[\tanh(0.5t),\sin(0.7t),\cos(0.5t)]^{\top},  $\\  $ d_{4}(\bar{{x}}_{4},t)=0.5\|0.02{x}_{41}\|[\cos(0.5t),\sin(0.7t),\tanh(0.5t)]^{\top}, $ \\
	$ d_{5}(\bar{{x}}_{5},t)=0.5\|0.04{x}_{51}\|[\tanh(0.5t),\sin(0.7t),\cos(0.5t)]^{\top}, $ \\  $ d_{6}(\bar{{x}}_{6},t)=0.5\|0.03{x}_{61}\|[\sin(0.5t),\cos(0.7t),\cos(0.5t)]^{\top}.  $
	
	{The control problem is to design a fixed-time controller such that the agents can form a desired formation in a dynamical environment with moving obstacles.
		Each agent can detect its surrounding environment using vehicle-mounted sensors with sensing range $10$m. In this simulation,} the initial positions of the agents are $x_{11}=[6,2,0]^{\top}$m, $x_{21}=[3,3+\sqrt{3},0]^{\top}$m, $ x_{31}=[-3,3+\sqrt{3},0]^{\top}$m, $x_{41}=[-7,-1,0]^{\top}$m, $ x_{51}=[-3,-3-\sqrt{3},0]^{\top}$m, $ x_{61}=[3,-3-\sqrt{3},0]^{\top}$m, and
	{the environmental obstacles have time-varying positions as} {$O_{1}=[0,2-\cos t,23]^{\top}$m, $ O_{2}=[1,-2-\cos t,28]^{\top}$m, $O_{3}=[-1.5,-\cos t,10]^{\top}$m, $ O_{4}=[1,-2-\cos t,7]^{\top}$m}. The radius of its repulsive zone in \eqref{2} is $d=2$m.
	
	The design parameters in the estimator \eqref{11} are selected as $ K_{1}=0.4 $, $ K_{2}=0.6 $, $K_{3}=1$, $ r_{3}=6 $, $ r_{4}=5 $, $ r_{5}=3 $, $ r_{6}=5 $. In the calculation of the sliding mode $ {\mathcal{S}}_{i} $ in \eqref{S}, we choose $ \varrho=100 $, $ c_{1}=2 $, $ c_{2}=0.2 $, $ \beta_{1}=\beta_{2}=0.6 $, $ \phi_{s}=0.01 $, $ r_{0}=0.9 $, $ r_{1}=1.2 $, $ r_{2}=0.6 $. The parameter $ \lambda_{f} $ in \eqref{v2} is chosen as $ \lambda_{f}=1 $.
	The setup for the RBF neural network used in \eqref{33} is as follows. Six neurons are contained in the neural network, and the sigmoid basis functions are applied with the center of the receptive field $ \mu_{ik}=k-3 $ and the width of the Gaussian function $ \psi_{ik}=\sqrt{2} $ for $ i=1,\dots,6 $, $ k=1,\dots,6 $, {where $ k $ denotes the $ k $th element of the basis function vector $ \phi_{i} $ and $ i $ {represents} the $ i $th agent.} Furthermore, we choose the design parameters $k_{1i}^{0}=k_{2i}^{0}=0.1$, $k_{1i}^{M}=k_{2i}^{M}=100$, $c_{s1}=c_{s2}=0.01$, $ \Gamma_{i}=I $, $ \gamma_{3i}=1 $ in the computation of ${u}_{i}$ in (\ref{u}), $ \hat W_{i} $ in \eqref{ad1} and $\hat\delta_{i} $ in \eqref{ad2}. Next, we give the initial value of the adaptive parameters $ \hat W_{i}(0)= 0$, $ \hat\delta_{i}(0)=0.1 $. 
	
	{To demonstrate the feasibility and flexibility of the proposed control scheme, two simulation scenarios are taken, where the first scenario considers
		the cooperative formation that maintains a prefixed relative distance between the virtual leader and each agent, and in the second scenario, we consider the cooperative formation problem using relative positions.}
	
	\subsection{Cooperative Tracking using Relative Distance}
	
\begin{figure}[!tp]\centering
		\includegraphics[width=8.2cm,trim=40 0 40 15, clip]{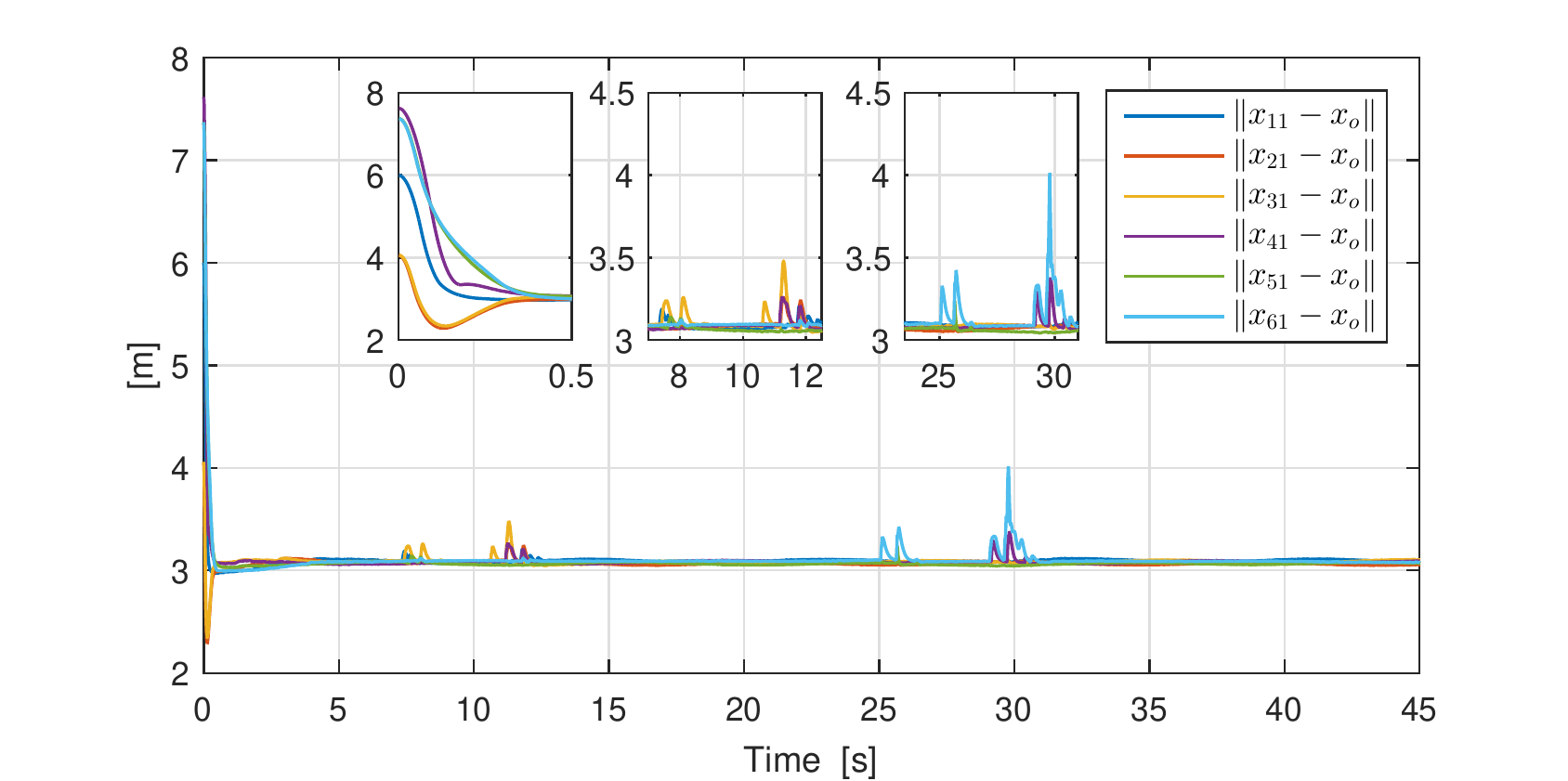}
		\caption{Responses of the distances between agents and leader.} 
		\label{fig:2}
	\end{figure}

\begin{figure}[!tp]\centering
		\includegraphics[width=8.2cm,trim=40 0 40 15, clip]{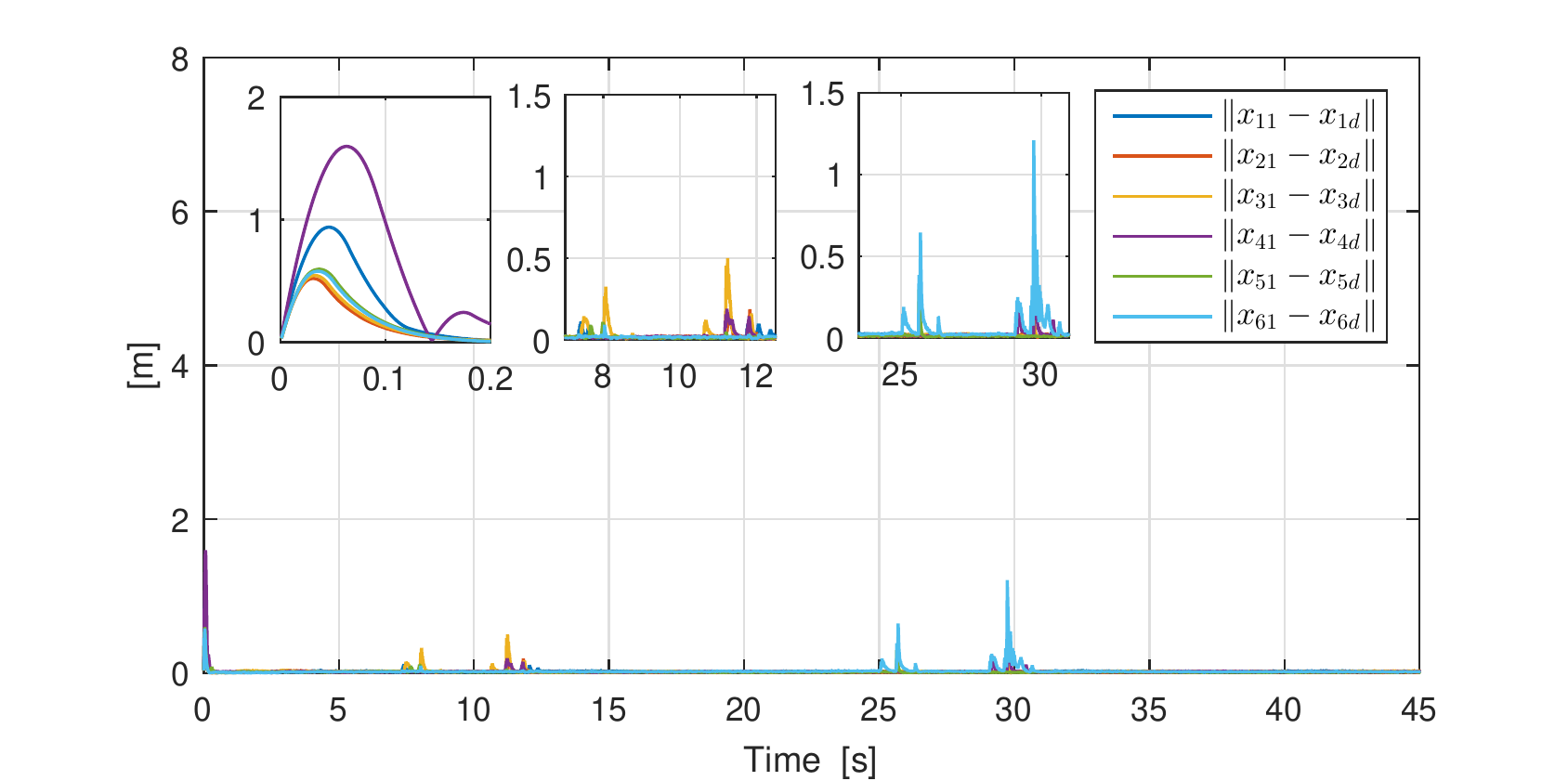}
		\caption{Responses of the six agents' tracking errors.} 
		\label{fig:3}
	\end{figure}

\begin{figure}[!tp]\centering
		\includegraphics[width=8.2cm,trim=40 0 40 15, clip]{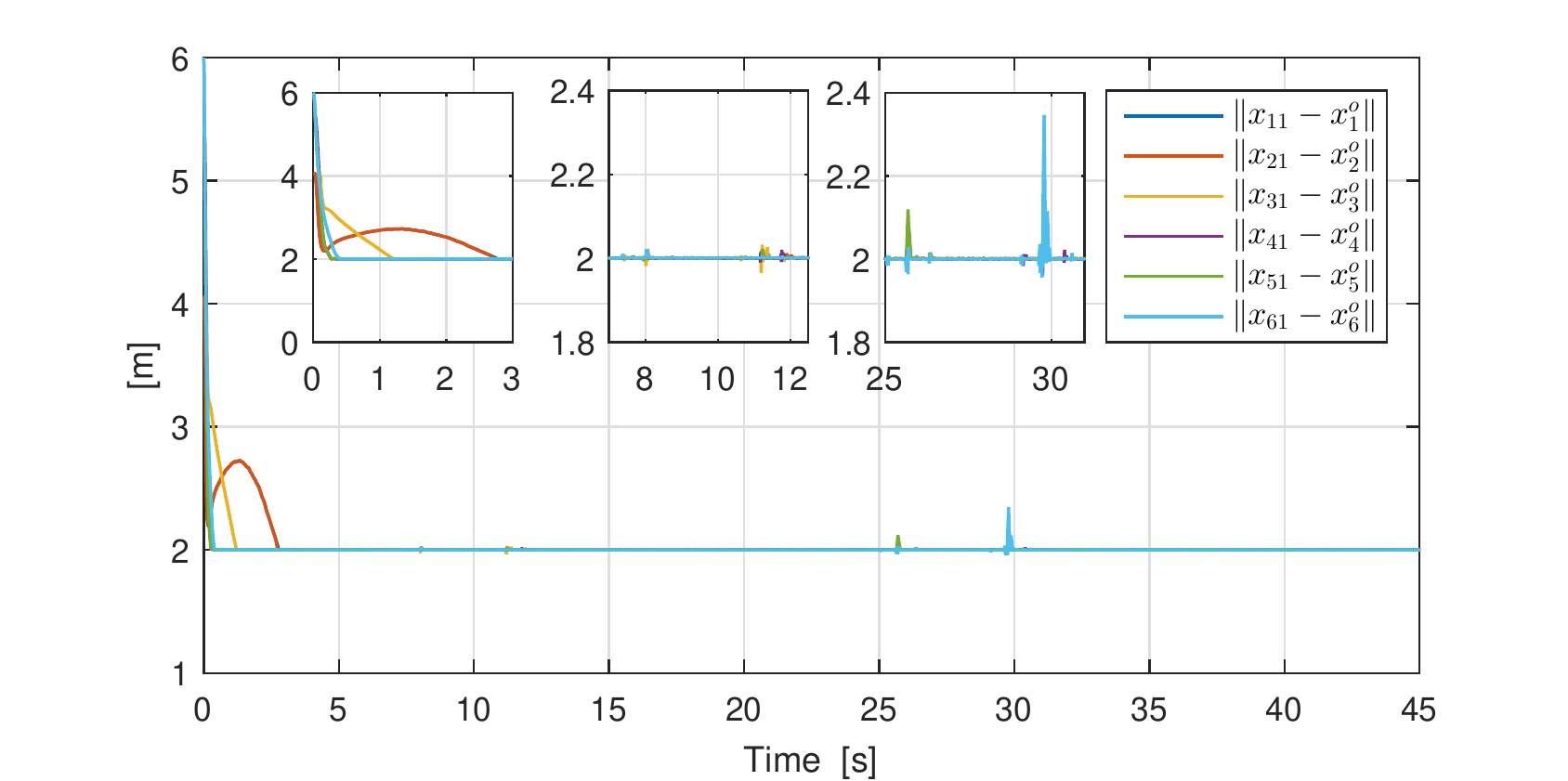}
		\caption{Responses of the distances between agents and obstacles.} 
		\label{fig:4}
	\end{figure}

\begin{figure} [!tp]\centering
		\includegraphics[width=8.2cm,trim=50 15 50 50, clip]{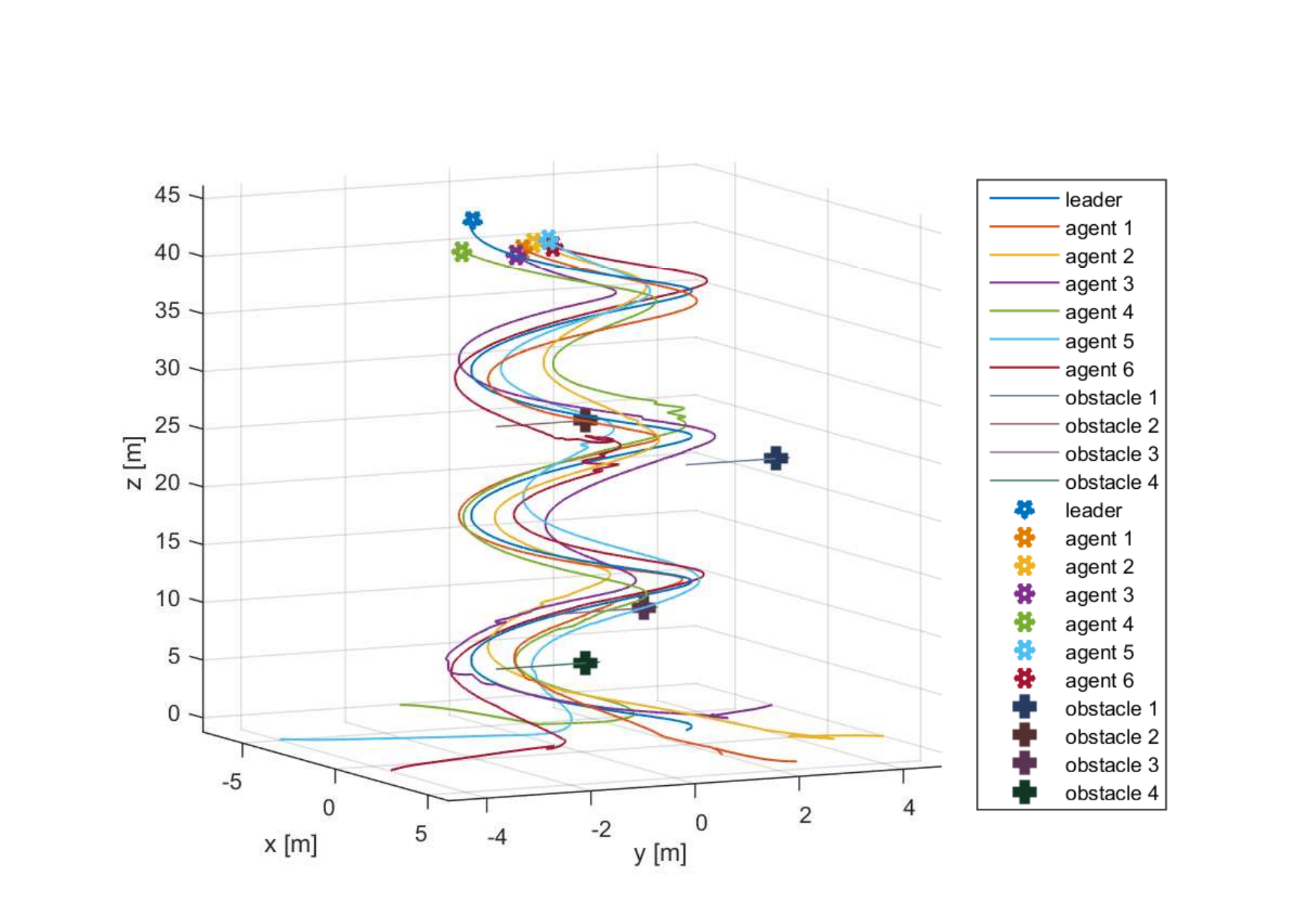}
		\caption{Trajectories of the six agents in the environment with four obstacles.} 
		\label{fig:5}
	\end{figure}

	{In this scenario, the simulation result is shown in Fig.~\ref{fig:2}, where the relative distances between six agents and the leader are presented in three-dimensional space. Observe that in the time intervals of {7s$\sim$8s, 11s$\sim $12s, 25s$\sim $26.5s, 29s$\sim $30.6s}, some of the relative distances in Fig.~\ref{fig:2} and  the tracking errors in Fig.~\ref{fig:3} deviate their original steady trajectories. This is because the corresponding agents encounter certain environmental obstacles, and their collision avoidance behaviors take place and have a higher priority in the desired velocity \eqref{v}. The controllers then steer these agents to move away from obstacles, resulting in a deviation from their desired trajectories for formation task.  
		
		From Fig. \ref{fig:4}  we can see that each agent can approximately maintain a distance $2$m away from the nearest object in its repulsive zone. The resulting formation is flexible as only the relative distance, rather than the relative position, from each agent to the virtual leader is fixed. The agents just move to certain positions on the hypersphere of the leader with a radius of $ d_{i0}=3 $m. The trajectories of six agents in the formation are presented in Fig.~\ref{fig:5}, which shows that all the agents can follow the leader with a distance $ d_{i0}=3 $m and adjust their motions to avoid collisions. Thus, the proposed control law in \eqref{u} works well in this cooperative formation task with the obstacle avoidance requirement.}
	

	\subsection{Cooperative Tracking using Relative Positions}
	
	{In this scenario, the formation task is performed based on relative position from the virtual leader to each agent. In the simulation, we set the expected relative positions as} \\$ x_{r11}=[\frac{3\sqrt{3}}{2},\frac{3}{2},0]^{\top} $m, $ x_{r21}=[0,3,0]^{\top} $m,\\ $ x_{r31}=[-\frac{3\sqrt{3}}{2},\frac{3}{2},0]^{\top} $m, $ x_{r41}=[-\frac{3\sqrt{3}}{2},-\frac{3}{2},0]^{\top} $m, \\$ x_{r21}=[0,-3,0]^{\top} $m, $ x_{r51}=[\frac{3\sqrt{3}}{2},-\frac{3}{2},0]^{\top} $m. 
	
\begin{figure}[!tp]\centering
		\includegraphics[width=8.2cm,trim=20 0 40 14, clip]{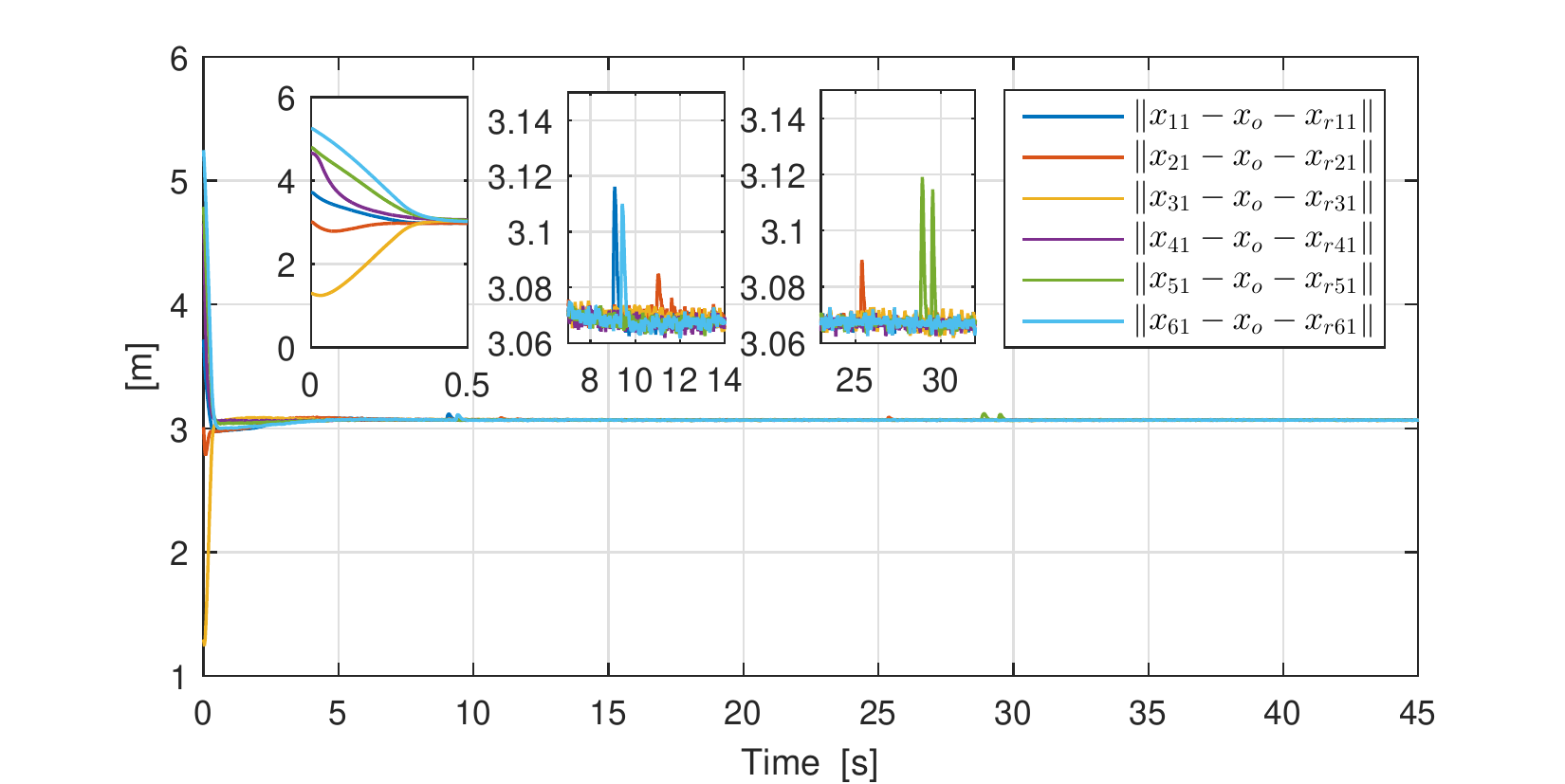}
		\caption{Responses of the distances between agents and leader.} 
		\label{fig:6}
	\end{figure}

\begin{figure} [!tp]\centering
		\includegraphics[width=8.2cm,trim=20 0 40 14, clip]{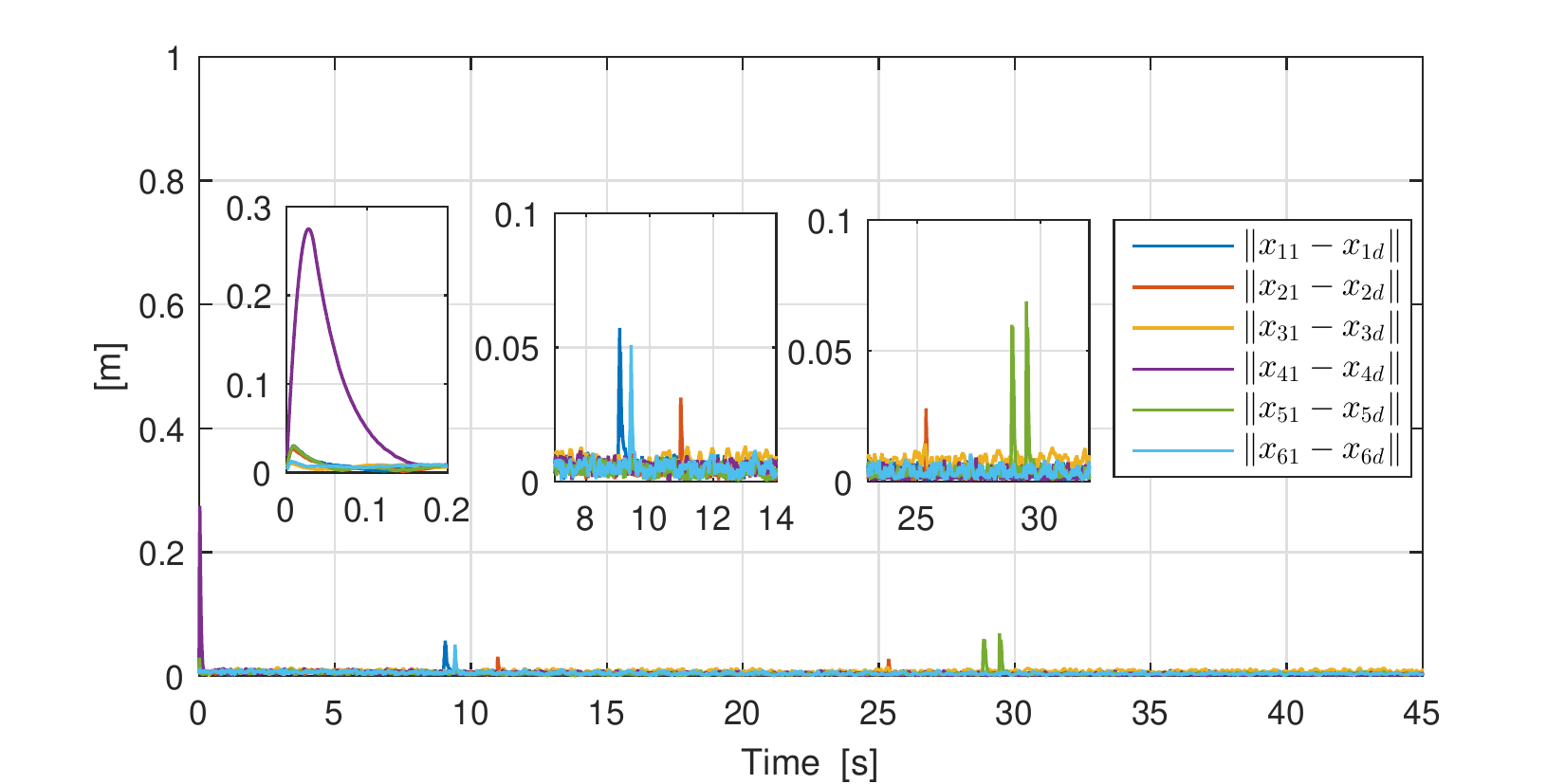}
		\caption{Responses of the six agents' tracking errors.} 
		\label{fig:7}
	\end{figure}

\begin{figure} [!tp]\centering
		\includegraphics[width=8.2cm,trim=20 0 40 14, clip]{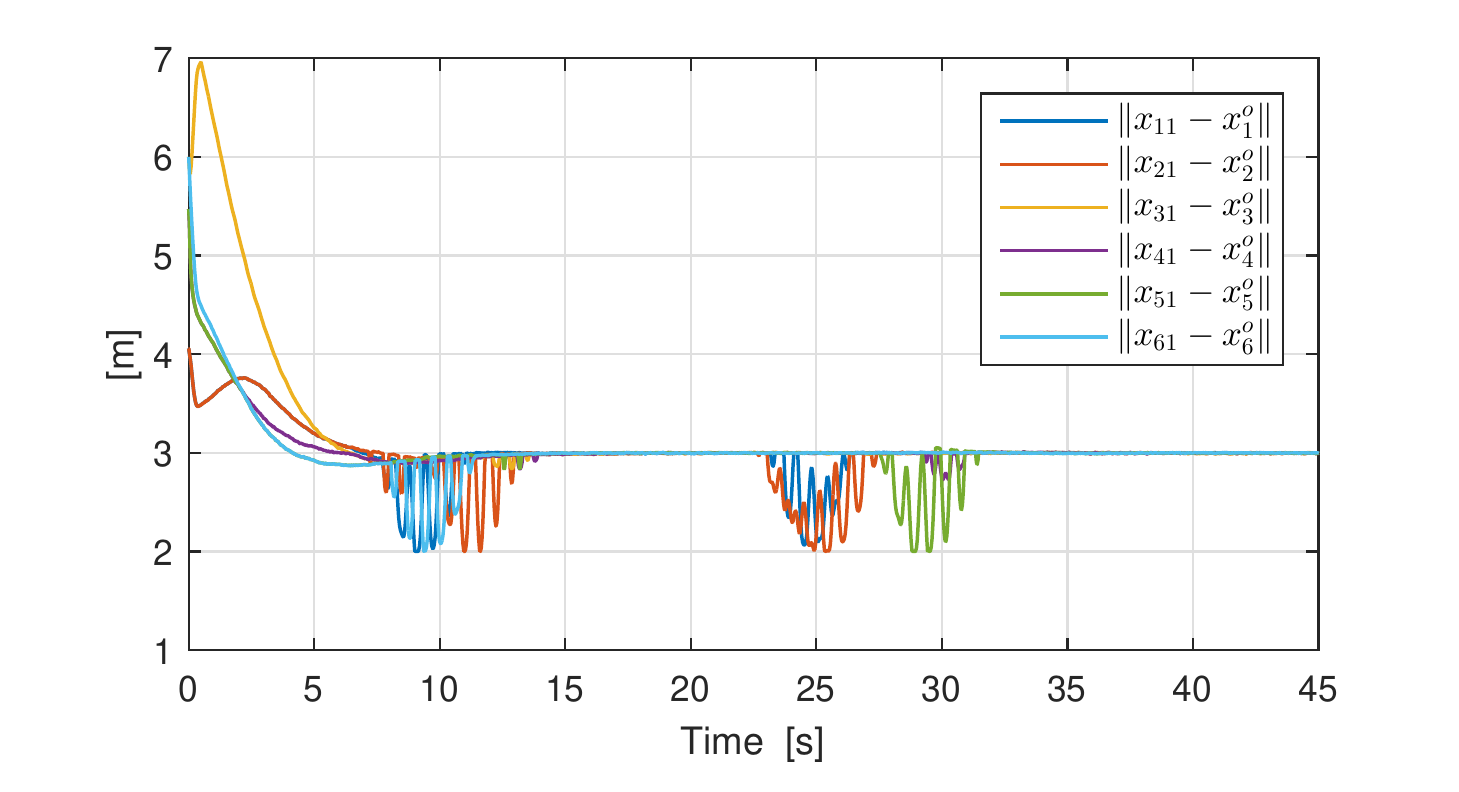}
		\caption{Responses of the distances between agents and obstacles.} 
		\label{fig:8}
	\end{figure}

\begin{figure} [!tp]\centering
		\includegraphics[width=8.2cm,trim=15 10 20 25, clip]{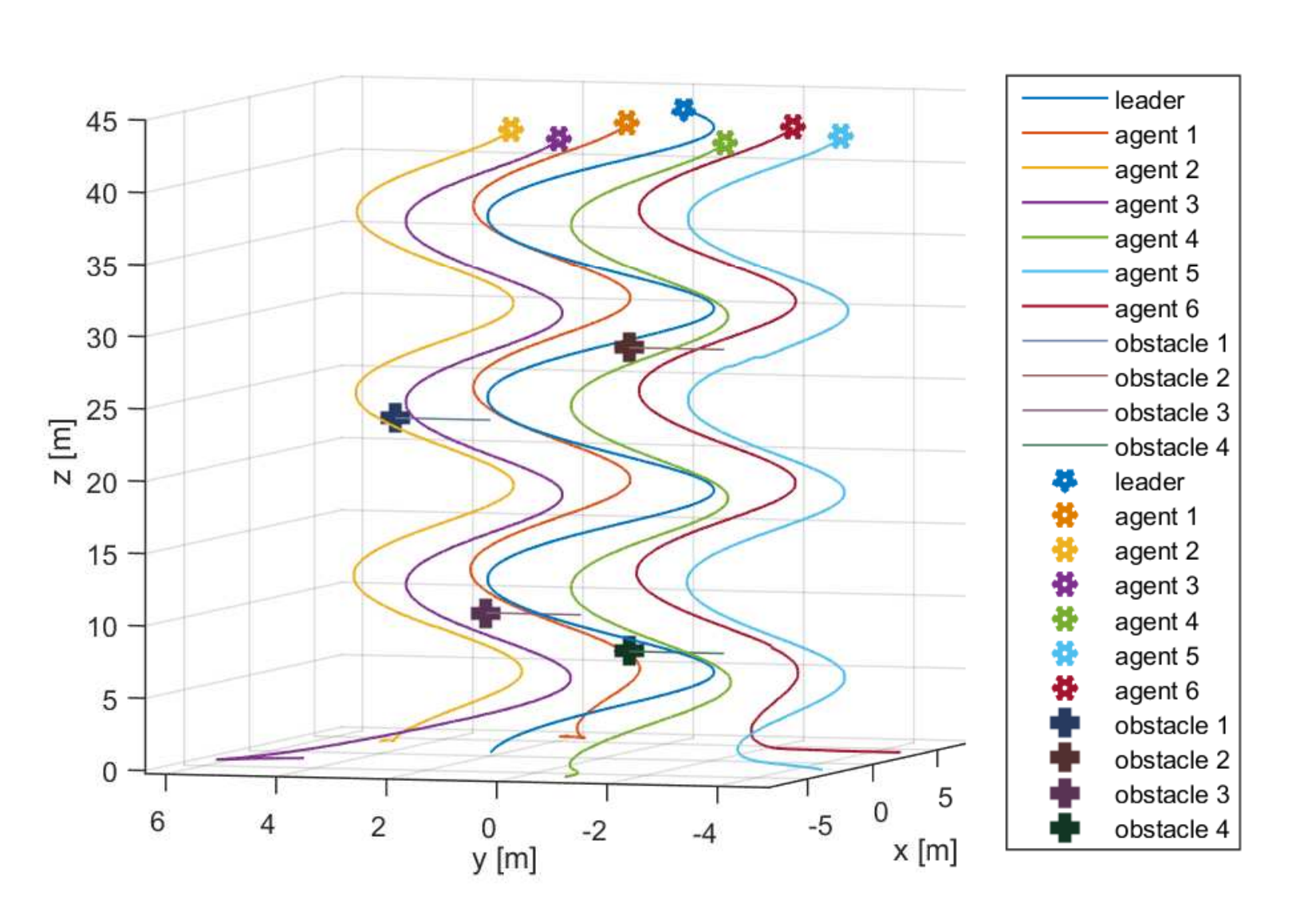}
		\caption{Trajectories of the six agents in the environment with four obstacles.}
		\label{fig:9}
	\end{figure}

	{The task error of each agent  in Fig. \ref{fig:6} shows the gap between the expected relative positions and the real relative positions at time $t$. In the time intervals of 9s$\sim$11.6s and 25s$\sim $30s, the obstacle avoidance behavior of the agents 1, 6, 2 and agents 2, 5 take place, respectively. The effects are also shown in the tracking error curves of the controlled agents in Fig.~\ref{fig:7}. Fig.~\ref{fig:8} shows the relative distances between each agent and its nearest object, which indicates that each agent   maintains a distance $ d=2 $m away from the environmental obstacles and the other agents. It is shown in Fig.~\ref{fig:9}  that all the agents follow the leader with a distance $ d_{i0}=3 $m and keep the expected relative positions. Therefore, the proposed algorithm \eqref{u} can {fulfil} the relative position based cooperative formation task in this dynamical environment.
	}
	
	{
		\subsection{Comparisons with Previous Works}
		To further show the performance of the proposed approach, we also simulate the methods in the previous works \cite{huang2019adaptive,zhou2017finite} and  make a comparison between them and our approach. In the simulations, we select the same initial conditions and design parameter values, and then compare the methods based on two indicators: the settling time  $ T_{s_{*}} $ and the overall precision index
		$
		I_{e}:=(\sum_{i=1}^{n}\|\tilde{x}_{i1}\|^{2})^{\frac{1}{2}},
		$
		where $ n $ denotes the number of agents. 
		The performances of different controllers \eqref{u} in this paper, (17) in \cite{huang2019adaptive} and (13) in \cite{zhou2017finite} are compared in Table~\ref{tab2}.
		\begin{table}[h]
			\centering\caption{Performance comparison}
			\begin{tabular}{c|c|c|c}
				Controller &Settling time $ T_{s_{*}} $ &$I_{e}$ at $ t=20s $ & $I_{e}$ at $ t=45s $\\ \hline 
				\eqref{u}&0.5s & $ 1.126\times 10^{-2} $&$ 8.523\times 10^{-3} $ \\ \hline 
				(17) in \cite{huang2019adaptive}&15s  &$ 1.015\times 10^{-1} $  & $ 1.027 \times 10^{-2} $\\ \hline 
				(13) in \cite{zhou2017finite}&8.5s &$ 2.053\times 10^{-2} $  &$ 1.442\times 10^{-2} $\\ 
			\end{tabular}\label{tab2}
		\end{table}
		
		From this table, we can see that the convergence rate of the proposed scheme is significantly faster than the other two. This is due to the strength of the fixed-time control strategy implemented in this work, compared to the finite-time control schemes in \cite{huang2019adaptive,zhou2017finite}.
		Meanwhile, the proposed method also achieves a higher control accuracy in terms of the index $I_e$ than the approaches in \cite{huang2019adaptive,zhou2017finite}. One important reason is that we take the advantage of the design parameter $\varrho$ in \eqref{S}, providing us an extra freedom to adjust the control precision.}
	
	\section{Conclusions}
	\label{sec:Conclusions}
	This paper investigated the problem of fixed-time cooperative control for a network of second-order nonlinear multi-agent systems. We developed a novel fixed-time behavioral control scheme for the multi-agent systems over undirected graphs. The proposed control protocol has been proven to overcome the disadvantage of the typical behavioral control, which works in a centralized manner and achieve the fixed-time convergence of tracking errors theoretically. The approximation of RBFNNs was adopted to solve the uncertainties of the system and the partial term of the sliding mode. Two behaviors for each autonomous agent were carefully defined and properly arranged in priority aiming at achieving the cooperative formation behavior. By employing these techniques, a fixed-time behavior controller has been designed for each agent such that all the agents can converge to the desired formation and avoid collisions, and the error signals achieved fixed-time convergence. Numerical simulation results have shown the effectiveness of the proposed approaches. The proposed control strategy will further be implemented using practical autonomous vehicles.

	\section*{Appendix}
	\renewcommand{\thesubsection}{\Alph{subsection}}

	\subsection{Proof of Lemma~\ref{l1}}\label{Ap:lem2}
	First, to examine the convergence property of $\tilde{\rho}_{io}$, we define a Lyapunov function candidate as follows 
	$$ V_{i,o}=\frac{1}{2}\gamma_{i,o}\tilde{\rho}_{io}^{2}, 
	$$
	where $\gamma_{i,o}>0$ is a design parameter. Differentiating
	$V_{i,o}$ with respect to time and using the {desired} velocity (\ref{v1}), it yields
	\begin{align*}
		\nm\dot{V}_{i,o}=&-\gamma_{i,o}\tilde{\rho}_{io}(J_{io}\dot{x}_{io}{+J_{i}^{o}\dot{x}_{i}^{o})}
		\\ \nm=&-\gamma_{i,o}\tilde{\rho}_{io}
		J_{io}J_{io}^{\dag}{[\lambda_{io}\alpha_{io}(\tilde{\rho}_{io})-J_{i}^{o}\dot{x}_{i}^{o}]-\gamma_{i,o}\tilde{\rho}_{io}J_{i}^{o}\dot{x}_{i}^{o}}\\
		=&-\gamma_{i,o}\lambda_{i,o}\tilde{\rho}_{io}\alpha_{io}(\tilde{\rho}_{io}).
	\end{align*}
	Then two cases are discussed according to the structure of \eqref{6}.
	
	\textbf{Case A}: If $\sigma_{1,i}=0$ or
	$\sigma_{1,i}\neq0$, $|\tilde{\rho}_{io}|>\phi_{s}$,
	where $\sigma_{1,i} (\dot{\tilde{\rho}}_{io},\tilde{\rho}_{io}):=\dot{\tilde{\rho}}_{io}+c_{0}(\beta_{1}\tilde{\rho}_{io}^{[r_{1}]}+\beta_{2}\tilde{\rho}_{io}^{[r_{2}]})^{[r_{0}]}$, we have
	\begin{align}\label{20}
		\dot{V}_{i,o} &=-{\gamma_{i,o}\lambda_{i,o}}\left[\beta_{1}\tilde{\rho}_{io}^{2(\frac{r_{1}r_{0}+1}{2r_{0}})}+\beta_{2}\tilde{\rho}_{io}^{2(\frac{r_{2}r_{0}+1}{2r_{0}})}\right]^{r_{0}} 
		\nm\\
		&=-\left[\gamma^{i}_{o1}V_{i,o}^{\frac{r_{1}r_{0}+1}{2r_{0}}}+\gamma^{i}_{o2}V_{i,o}^{\frac{r_{2}r_{0}+1}{2r_{0}}}\right]^{r_{0}},
	\end{align}
	where  
	\begin{align*}
		\gamma^{i}_{o1}:= ({\gamma_{i,o}\lambda_{i,o}})^{\frac{1}{r_{0}}}
		\left( \frac{2}{\gamma_{i,o}}\right)^{\frac{r_{1}r_{0}+1}{2r_{0}}}\beta_{1}, 
		\\
		\gamma^{i}_{o2}:=({\gamma_{i,o}\lambda_{i,o}})^{\frac{1}{r_{0}}}
		\left(\frac{2}{\gamma_{i,o}}\right)^{\frac{r_{2}r_{0}+1}{2r_{0}}}\beta_{2}.
	\end{align*}  
	Then, it follows from Lemma~\ref{ld1} that $\tilde{\rho}_{io}$
	converges to 0 in a fixed time $T_{i,o}$ with
	$$ T_{
		i,o}\leq\frac{1}{(\gamma^{i}_{o1})^{r_{0}}(r_{1}r_{0}-1)}+\frac{1}{(\gamma^{i}_{o2})^{r_{0}}(1-r_{2}r_{0})}, \forall \ \tilde{\rho}_{io}(0)\in\bR^{3}. $$
	
	\textbf{Case B}: If $\sigma_{1,i}\neq0$ and
	$|\tilde{\rho}_{io}|\le\phi_{s}$, we obtain
	\begin{align} \label{eq:dotVio}
		\dot{V}_{i,o}
		&=-\gamma_{i,o}\lambda_{i,o}(\wp_{1}\tilde{\rho}_{io}^{2}+\wp_{2}|\tilde{\rho}_{io}|^{3})\nm\\
		&\le -\gamma_{i,o}\lambda_{i,o}\wp_{1}\tilde{\rho}_{io}^{2}+\gamma_{i,o}\lambda_{i,o}|\wp_{2}|\phi_{s}^{3}\nm\\&
		=-\gamma^{i}_{o3}V_{i,o}+\gamma^{i}_{o4}.
	\end{align}
	where the definition of
	$\alpha_{io}(\tilde{\rho}_{io})$ in \eqref{6} is used, and $ \gamma^{i}_{o3}:=2\lambda_{i,o}\wp_{1} $, $ \gamma^{i}_{o4}:= \gamma_{i,o}\lambda_{i,o}|\wp_{2}|\phi_{s}^{3}$. From \eqref{eq:dotVio}, we conclude that for any $ {\rho}_{io}(0) \in \bR_{\ge 0}$, there exists a constant $ \varsigma_{\rho i} \in \bR_{> 0}$, related to $ {\rho}_{io}(0) $, such that $|\tilde{\rho}_{io}|\le\varsigma_{\rho i} $. Furthermore, it implies from $|\tilde{\rho}_{io}|\leq\phi_{s}$ 
	that $\tilde{\rho}_{io}$ can converge to the region $|\tilde{\rho}_{io}|\leq\phi_{s}$ in a fixed time.
	
	We combine the analysis in both cases and conclude that
	$\tilde{\rho}_{io}$ converges to the region $|\tilde{\rho}_{io}|\leq\phi_{s}$ 
	in a fixed time using {desired} velocity \eqref{v1}, meaning that $ \|x_{i1}-x_{i}^{o}\|\ge \sqrt{d^{2}-2\phi_{s}} $.
	
	\subsection{Proof of Lemma~\ref{l8}} \label{Ap:lem8}
	Denote $H: = L + B$, which is symmetric and positive definite \cite{godsil2013algebraic}. We
	consider a Lyapunov candidate
	\begin{align}\label{12}
		{V_{e}(\bar{x}_{1}):=\frac{1}{2}\bar{x}_{1}^{\top} (H\otimes I_{3}) \bar{x}_{1},}
	\end{align}
	where $\bar{x}_{1}(t) =[\bar{x}_{11}^{\top}(t),\dots,\bar{x}_{n1}^{\top}(t)]^{\top} $ {with $ \bar{x}_{i1}$ defined in \eqref{eq:bar_x}. Denote $\hat{x}_{1}(t) =[\hat{x}_{11}^{\top}(t),\dots,\hat{x}_{n1}^{\top}(t)]^{\top} $. } The time derivative of $ V_{e} $ is computed along \eqref{11} as{\begin{align}\label{13}
			\dot{V}_{e}=&\bar{x}_{1}^{\top}{(H\otimes I_{3})}(\dot{\hat{x}}_{1}-\1_{n}\otimes \dot{x}_{o})
			\nm\\
			=&{\bar{x}_{1}^{\top} (H\otimes I_{3})[-K_{1}((H\otimes I_{3})\bar{x}_{1})^{[\frac{r_{3}}{r_{4}}]}-\1_{n}\otimes \dot{x}_{o}}
			\nm\\ & {-K_{2}((H\otimes I_{3})\bar{x}_{1})^{[\frac{r_{5}}{r_{6}}]}-K_{3}{\rm{sgn}}((H\otimes I_{3})\bar{x}_{1})] }
			\nm\\
			\le&{-K_{1}\bar{x}_{1}^{\top}(H\otimes I_{3})((H\otimes I_{3})\bar{x}_{1})^{[\frac{r_{3}}{r_{4}}]}}
			\nm\\& {-K_{2}\bar{x}_{1}^{\top}(H\otimes I_{3})((H\otimes I_{3})\bar{x}_{1})^{[\frac{r_{5}}{r_{6}}]}} \nm\\&{-{(K_{3}-{\sup_{t\ge 0}\|\dot{x}_{o}(t)\|_{\infty}})}\|(H\otimes I_{3})\bar{x}_{1}\|_{1}.}
		\end{align}Then the bound in \eqref{a3.1} is applied for this proof.
		Using Lemma 6 for \eqref{13}, it becomes\begin{align}\label{13a}
			\dot{V}_{e}\le&{-K_{1}{(3N)^{1-\frac{r_{3}+r_{4}}{2r_{4}}}}(\bar{x}_{1}^{\top}(H\otimes I_{3})(H\otimes I_{3})\bar{x}_{1})^{\frac{r_{3}+r_{4}}{2r_{4}}}} \nm\\&
			{-K_{2}(\bar{x}_{1}^{\top}(H\otimes I_{3})(H\otimes I_{3})\bar{x}_{1})^{\frac{r_{5}+r_{6}}{2r_{6}}}}.
		\end{align}
		According to the properties of symmetrical positive definite matrix, we have $ HX= \lambda_{H}X$, where $ \lambda_{H}:={\rm diag}\{\lambda_{H1}, \dots, \lambda_{Hn}\}>0 $, $ X:=[X_{1},\dots,X_{n}]\in\bR^{n\times n} $, $ X_{1},\dots,X_{n} \in\bR^{n}$ are the eigenvectors of $ H $ corresponding to its eigenvalues $ \lambda_{H1}, \dots, \lambda_{Hn}$, which can be chosen as a set of orthogonal bases of $ \bR^{n} $ and satisfies $ X^{\top}X=XX^{\top}=I_{n} $. Based on the above property, we can rewrite \eqref{13a} as follows:\begin{align}
			\dot{V}_{e}
			\le&- \tilde{K}_{1}{V}_{e}^{\tilde r_{1}}- \tilde{K}_{2}{V}_{e}^{\tilde r_{2}},
		\end{align}
		where $\tilde{K}_{1}:= {K}_{1}{(3N)^{1-\frac{r_{3}+r_{4}}{2r_{4}}}}(2\frac{\lambda_{\min}^{2}(H)}{\lambda_{\max}(H)})^{\frac{1}{\tilde r_{1}}} $, $ \tilde{K}_{2}:= {K}_{2}(2\frac{\lambda_{\min}^{2}(H)}{\lambda_{\max}(H)})^{\frac{1}{\tilde r_{2}}}$, $\tilde r_{1}:= \frac{r_{3}+r_{4}}{2r_{4}}$ and $\tilde r_{2}:=\frac{r_{5}+r_{6}}{2r_{6}} $.}
	
	Then, it follows from Lemma~\ref{ld1} that for any $ \hat{x}_{i1}(0)$ and $x_{o}(0)) $, there exists 
	a convergence time 
	{\begin{equation*} 
			T_{e}:=\frac{1}{ \tilde{K}_{1} (\tilde r_{1}-1)}+\frac{1}{\tilde{K}_{2}(1-\tilde r_{2})},
	\end{equation*}}
	such that $\hat{x}_{i1}(t)\equiv x_{o}(t) $ when $ t\ge T_{e} $.
	
	\subsection{Proof of Lemma~\ref{l4}} \label{Ap:lem4}
	This lemma can be proved using the similar procedure as the proof of Lemma~\ref{l1}. However, we use a different Lyapunov function candidate
	$$ V_{f}:=\frac{1}{2}\gamma_{f}\tilde{\rho}_{f}^{\top}\tilde{\rho}_{f},\ V_{i,f}:=\frac{1}{2}\gamma_{f}\tilde{\rho}_{if}^{2}, $$
	where $\gamma_{f}>0$ is the design parameter in \eqref{v2}. Taking the derivative of
	$V_{f}$ with respect to time then leads to
	\begin{align*}
		\nm\dot{V}_{f}=&-\gamma_{f}\tilde{\rho}_{f}^{\top}\dot{\tilde{\rho}}_{f}\nm\\=&-\gamma_{f}\tilde{\rho}_{f}^{\top}(J_{f}\dot{x}_{f}+J_{\hat{f}}\dot{\hat{x}}_{1})
		\\ \nm=&-\gamma_{f}\tilde{\rho}_{f}^{\top}
		J_{f}J^{\dagger}_{f} [\Lambda_{f}\alpha_{f}(\tilde{\rho}_{f})-J_{\hat{f}}\dot{\hat{x}}_{1}]-\gamma_{f}\tilde{\rho}_{f}^{\top}J_{\hat{f}}\dot{\hat{x}}_{1}\\
		=&-\gamma_{f}\lambda_{f}\tilde{\rho}_{f}^{\top}\alpha_{f}(\tilde{\rho}_{f})\\
		=&-\sum_{i=1}^{n}\gamma_{f}\lambda_{f}\tilde{\rho}_{if}\alpha_{if}(\tilde{\rho}_{if}).
	\end{align*}
	For each individual agent $i$, we have
	$
	\nm\dot{V}_{i,f}=-\gamma_{f}\lambda_{f}\tilde{\rho}_{if}\alpha_{if}(\tilde{\rho}_{if}).
	$
	Then the rest of this proof follows similarly as the proof of Lemma~\ref{l1}, and we omit the details due to the limited space.
	
	\subsection{Proof of Theorem~\ref{taskl4}}\label{Ap:task14}
	{According to the definitions of $\alpha_{io}(\tilde{\rho}_{io})$ and $\alpha_{if}(\tilde{\rho}_{if})$ in \eqref{6} and \eqref{if}, four cases are discussed.}
	
	\textbf{Case A}: $ \sigma_{1,i} (\dot{\tilde{\rho}}_{io},\tilde{\rho}_{io})\neq0$, $ |\tilde{\rho}_{io}(t)| \leq \phi_{s} $ and $ \sigma_{1,i} (\dot{\tilde{\rho}}_{if},\tilde{\rho}_{if})\neq0$, $ |\tilde{\rho}_{if}(t)| \leq \phi_{s} $. It is immediate that $ \|x_{i1}-x_{i}^{o}\|\ge \sqrt{d^{2}-2\phi_{s}} $ and {$ \|x_{i1}-\hat{x}_{i1}\|\ge \sqrt{d_{i0}^{2}-2\phi_{s}} $} for all $i \in V$ and $t\geq 0$.
	
	\textbf{Case B}: $ \sigma_{1,i} (\dot{\tilde{\rho}}_{io},\tilde{\rho}_{io})=0$ or $ \sigma_{1,i} (\dot{\tilde{\rho}}_{io},\tilde{\rho}_{io})\neq0$, $ |\tilde{\rho}_{io}(t)| > \phi_{s} $ and $ \sigma_{1,i} (\dot{\tilde{\rho}}_{if},\tilde{\rho}_{if})\neq0$, $ |\tilde{\rho}_{if}(t)| \leq \phi_{s} $. According to Lemma \ref{l1}, there exists a settling time $T_{i,o}>0$ such that $ \|x_{i1}-x_{i}^{o}\|\ge \sqrt{d^{2}-2\phi_{s}} $ and {$ \|x_{i1}-\hat{x}_{i1}\|\ge \sqrt{d_{i0}^{2}-2\phi_{s}} $} for all $i\in V$ and $t\geq T_{i,o}$ .
	
	\textbf{Case C}: $ \sigma_{1,i} (\dot{\tilde{\rho}}_{io},\tilde{\rho}_{io})\neq0$, $ |\tilde{\rho}_{io}(t)| \leq \phi_{s} $ and $ \sigma_{1,i} (\dot{\tilde{\rho}}_{if},\tilde{\rho}_{if})=0$ or $ \sigma_{1,i} (\dot{\tilde{\rho}}_{if},\tilde{\rho}_{if})\neq0$, $ |\tilde{\rho}_{if}(t)| > \phi_{s} $. According to Lemma \ref{l4}, there exists a settling time $T_{i,f}>0$ such that $ \|x_{i1}-x_{i}^{o}\|\ge \sqrt{d^{2}-2\phi_{s}} $ and {$ \|x_{i1}-\hat{x}_{i1}\|\ge \sqrt{d_{i0}^{2}-2\phi_{s}} $} for all $i\in V$ and $t\geq T_{i,f}$.
	
	\textbf{Case D}: $ \sigma_{1,i} (\dot{\tilde{\rho}}_{io},\tilde{\rho}_{io})=0$ or $ \sigma_{1,i} (\dot{\tilde{\rho}}_{io},\tilde{\rho}_{io})\neq0$, $ |\tilde{\rho}_{io}(t)| > \phi_{s} $ and $ \sigma_{1,i} (\dot{\tilde{\rho}}_{if},\tilde{\rho}_{if})=0$ or $ \sigma_{1,i} (\dot{\tilde{\rho}}_{if},\tilde{\rho}_{if})\neq0$, $ |\tilde{\rho}_{if}(t)| > \phi_{s} $. 
	
	In the following proof, we focus on the analysis of this case. We design a Lyapunov function for agent $ i $ as 
	\begin{align}\label{taskl4_1}
	&V_{i,M}(\tilde{\rho}_{io},\tilde{\rho}_{if}):=V_{i,o}(\tilde{\rho}_{io})+V_{i,f}(\tilde{\rho}_{if}),\\ &
	V_{i,o}(\tilde{\rho}_{io}):=\frac{1}{2}\gamma_{i,o}\tilde{\rho}_{io}^{2}, \ V_{i,f}(\tilde{\rho}_{if}):=\frac{1}{2}\gamma_{f}\tilde{\rho}_{if}^{2},\nm
	\end{align}where $ \gamma_{i,o}$, $\gamma_{f}>0 $, $ \gamma_{i,o}\ge \frac{d_{i0}\gamma_{f}L_{0}}{d\phi_{s}}+ \gamma_{i,\varepsilon}$ with ${0<} \phi_{s}<\max_{i}\{{|\tilde{\rho}_{io}(0)|,|\tilde{\rho}_{if}(0)|}\}\le L_{0} $, $\gamma_{i,\varepsilon}>0$. $ \tilde{\rho}_{io}(0),\tilde{\rho}_{if}(0) $ are the initial values of $ \tilde{\rho}_{io}(t) $ and $\tilde{\rho}_{if}(t) $.
	
	When $\|\dot{x}_{io}+\dot{x}_{if}\|\neq 0 $, taking the time derivative of $ V_{i,M} $ along the desired velocity $\dot{x}_{id}$ in \eqref{v} gives
	\begin{align}\label{16}
	\dot{V}_{i,M}
	=&-\gamma_{i,o}\tilde{\rho}_{io}J_{io}[\dot{x}_{io}+(I-J_{io}^{\dag}J_{io})\dot{x}_{if}]{-\gamma_{i,o}\tilde{\rho}_{io}J_{i}^{o}\dot{x}_{i}^{o}}\nm\\
	&-\gamma_{f}\tilde{\rho}_{if}J_{if}[\dot{x}_{io}+(I-J_{io}^{\dag}J_{io})\dot{x}_{if}]-\gamma_{f}\tilde{\rho}_{if}J_{i\hat{f}}\dot{\hat{x}}_{i1} \nm\\
	=&-\gamma_{i,o}\tilde{\rho}_{io}J_{io}\{J_{io}^{\dag} \lambda_{io}(\beta_{1}\tilde{\rho}_{io}^{[r_{1}]}+\beta_{2}\tilde{\rho}_{io}^{[r_{2}]})^{[r_{0}]}\nm\\&+(I-J_{io}^{\dag}J_{io})J^{\dagger}_{if} [\lambda_{f}(\beta_{1}\tilde{\rho}_{if}^{[r_{1}]}+\beta_{2}\tilde{\rho}_{if}^{[r_{2}]})^{[r_{0}]}\nm\\&- J_{i\hat{f}}\dot{\hat{x}}_{i1}]\}
	-\gamma_{f}\tilde{\rho}_{if}J_{if}\{J_{io}^{\dag} [\lambda_{io}(\beta_{1}\tilde{\rho}_{io}^{[r_{1}]}+\beta_{2}\tilde{\rho}_{io}^{[r_{2}]})^{[r_{0}]}\nm\\&{-J_{i}^{o}\dot{x}_{i}^{o}]}+(I-J_{io}^{\dag}J_{io})J^{\dagger}_{if} [\lambda_{f}(\beta_{1}\tilde{\rho}_{if}^{[r_{1}]}+\beta_{2}\tilde{\rho}_{if}^{[r_{2}]})^{[r_{0}]}\nm\\&- J_{i\hat{f}}\dot{\hat{x}}_{i1}]\}-\gamma_{f}\tilde{\rho}_{if}J_{i\hat{f}}\dot{\hat{x}}_{i1},
	\end{align}
	where the definitions of $\alpha_{io}(\tilde{\rho}_{io})$ and $\alpha_{if}(\tilde{\rho}_{if})$ in \eqref{6} and \eqref{if} are used. 
	We can rewrite \eqref{16} as
	\begin{align}\label{18}
	\dot{V}_{i,M}\le&-\gamma_{i,o}\lambda_{io}|\tilde{\rho}_{io}||(\beta_{1}\tilde{\rho}_{io}^{[r_{1}]}+\beta_{2}\tilde{\rho}_{io}^{[r_{2}]})^{[r_{0}]}|\nm\\&+\gamma_{i,o}\lambda_{f}|J_{io}J^{\dagger}_{if}||\tilde{\rho}_{io}||(\beta_{1}\tilde{\rho}_{if}^{[r_{1}]}+\beta_{2}\tilde{\rho}_{if}^{[r_{2}]})^{[r_{0}]}|\nm\\&+\gamma_{f}\lambda_{io}|J_{if}J^{\dagger}_{io}||\tilde{\rho}_{if}||(\beta_{1}\tilde{\rho}_{io}^{[r_{1}]}+\beta_{2}\tilde{\rho}_{io}^{[r_{2}]})^{[r_{0}]}|\nm\\
	&+\gamma_{i,o}|\tilde{\rho}_{io}|\|J_{io}\|\|\dot{\hat{x}}_{i1}\|{+\gamma_{i,o}|\tilde{\rho}_{io}|\|J_{io}\|\|\dot{\hat{x}}_{i1}\|}\nm\\&+\gamma_{f}|\tilde{\rho}_{if}|\|J_{i\hat{f}}\|\|\dot{\hat{x}}_{i1}\|.
	\end{align}
	From this point, the proof goes in the following two directions. 
	\\
	
	\textbf{(a)} If $ |\tilde{\rho}_{io}|\ge |\tilde{\rho}_{if}| > \phi_{s}$, we have  $$|(\beta_{1}\tilde{\rho}_{io}^{[r_{1}]}+\beta_{2}\tilde{\rho}_{io}^{[r_{2}]})^{[r_{0}]}|\ge |(\beta_{1}\tilde{\rho}_{if}^{[r_{1}]}+\beta_{2}\tilde{\rho}_{if}^{[r_{2}]})^{[r_{0}]}|> \phi_{s}^{*},$$
	where $ \phi_{s}^{*}:= (\beta_{1}\phi_{s}^{r_{1}}+\beta_{2}\phi_{s}^{r_{2}})^{r_{0}}$.
	
	First, we prove that $ \tilde{\rho}_{io}(t) $ is bounded for an arbitrary bounded initial value $ \tilde{\rho}_{if}(0)$. Taking the time derivte of ${V}_{i,M} $ at $ t=0 $ and applying $ \gamma_{i,o}\ge \frac{\gamma_{f}L_{0}\sqrt{d_{i0}^{2}+2L_{0}}}{\phi_{s}\sqrt{d^{2}-2\phi_{s}}}+ \gamma_{i,\varepsilon}$ and $ 1<\frac{|\tilde{\rho}_{io}(t)|}{|\tilde{\rho}_{if}(t)|}\le\frac{ L_{0}}{\phi_{s}} $,   \eqref{18} becomes
	\begin{align}\label{d22}
	&\dot{V}_{i,M}(0)
	\le-\gamma_{i,o}\lambda_{i,o}|\tilde{\rho}_{io}(0)||(\beta_{1}\tilde{\rho}_{io}^{[r_{1}]}(0)+\beta_{2}\tilde{\rho}_{io}^{[r_{2}]}(0))^{[r_{0}]}|\nm\\&+\gamma_{i,o}\lambda_{f}L_{iof}|\tilde{\rho}_{io}(0)||(\beta_{1}\tilde{\rho}_{io}^{[r_{1}]}(0)+\beta_{2}\tilde{\rho}_{io}^{[r_{2}]}(0))^{[r_{0}]}|\nm\\&+\gamma_{f}\lambda_{io}L_{ifo}|\tilde{\rho}_{io}(0)||(\beta_{1}\tilde{\rho}_{io}^{[r_{1}]}(0)+\beta_{2}\tilde{\rho}_{io}^{[r_{2}]}(0))^{[r_{0}]}|\nm\\ 
	&+\frac{\gamma_{i,o}L_{io}L_{i\hat{x}}}{\phi_{s}^{*}}|\tilde{\rho}_{io}(0)||(\beta_{1}\tilde{\rho}_{io}^{[r_{1}]}(0)+\beta_{2}\tilde{\rho}_{io}^{[r_{2}]}(0))^{[r_{0}]}|\nm\\&{+\frac{\gamma_{f}L_{if}L_{i}^{o}}{\phi_{s}^{*}}|\tilde{\rho}_{if}(0)||(\beta_{1}\tilde{\rho}_{io}^{[r_{1}]}(0)+\beta_{2}\tilde{\rho}_{io}^{[r_{2}]}(0))^{[r_{0}]}|}\nm\\&+\frac{\gamma_{f}L_{if}L_{i\hat{x}}}{\phi_{s}^{*}}|\tilde{\rho}_{if}(0)||(\beta_{1}\tilde{\rho}_{io}^{[r_{1}]}(0)+\beta_{2}\tilde{\rho}_{io}^{[r_{2}]}(0))^{[r_{0}]}|\nm	\\
	&\le-\gamma_{i,o}\lambda_{iv}|\tilde{\rho}_{io}(0)||(\beta_{1}\tilde{\rho}_{io}^{[r_{1}]}(0)+\beta_{2}\tilde{\rho}_{io}^{[r_{2}]}(0))^{[r_{0}]}|\le 0,
	\end{align}
	\if0 where $ |J_{io}J^{\dagger}_{if}|\le L_{iof}$, $ |J_{if}J^{\dagger}_{io}| \le L_{ifo}$, $ \|J_{io}\|\le \sqrt{d^{2}+2L_{0}} $ and $ \|J_{i\hat{f}}\|\le \sqrt{d_{i0}^{2}+2L_{0}} $ are utilized according to Eqs. \eqref{jio} and \eqref{jif}, \fi
	where $ \lambda_{i,o} $ is designed to satisfy
	\begin{align}\label{lambdaio1}
	\lambda_{i,o}\ge\lambda_{i,o1}:=& \frac{\gamma_{i,o}\lambda_{f}L_{iof}}{\gamma_{i,\varepsilon}}+ \frac{\gamma_{i,o}L_{io}L_{i\hat{x}}}{\phi_{s}^{*}\gamma_{i,\varepsilon}}\nm\\&{+ \frac{\gamma_{f}L_{if}(L_{i\hat{x}}+L_{i}^{o})}{\phi_{s}^{*}\gamma_{i,\varepsilon}}} +\frac{\gamma_{i,o}{\lambda_{iv}}}{\gamma_{i,\varepsilon}},
	\end{align} {$ \lambda_{iv}>0 $} is an auxiliary design parameter, $ L_{iof}:= \sqrt{\frac{d^{2}+2L_{0}}{d_{i0}^{2}-2\phi_{s}}} $, $ L_{ifo}:= \sqrt{\frac{d_{i0}^{2}+2L_{0}}{d^{2}-2\phi_{s}}} $, $ L_{io}:=\sqrt{d^{2}+2L_{0}} $, {$ L_{i}^{o}:=\sup_{t\ge 0}\|\dot{x}_{i}^{o}(t)\|$,} $ L_{if}:=\sqrt{d_{i0}^{2}+2L_{0}} $ and 
	\begin{align}
	L_{i\hat{x}}:=& K_{1}\|\bar{\eta}_{i1}(0)\|+K_{2}\|\bar{\eta}_{i2}(0)\|+\sqrt{3}K_{3},\nm\\
	\bar{\eta}_{i1}(0)=& \sum_{j\in \mathcal{N}_{i}}a_{ij}(\bar{x}_{i1}^{[\frac{r_{3}}{r_{4}}]}(0)+\bar{x}_{j1}^{[\frac{r_{3}}{r_{4}}]}(0))+b_{i}\bar{x}_{i1}^{[\frac{r_{3}}{r_{4}}]}(0),\nm\\
	\bar{\eta}_{i2}(0)=& \sum_{j\in \mathcal{N}_{i}}a_{ij}(\bar{x}_{i1}^{[\frac{r_{5}}{r_{6}}]}(0)+\bar{x}_{j1}^{[\frac{r_{5}}{r_{6}}]}(0))+b_{i}\bar{x}_{i1}^{[\frac{r_{5}}{r_{6}}]}(0).\nm
	\end{align}
	The constraint $ \lambda_{i,o}\ge \lambda_{i,o1} $ is used to ensure $ \dot{V}_{i,M}(0)\le 0 $. From \eqref{d22}, we have 
	$ |\tilde{\rho}_{io}(\Delta t)|\le |\tilde{\rho}_{io}(0)|\le L_{0}$, for any   $ \Delta t \in \mathbb{R}_{>0}$, and similarly, we verify $ \dot{V}_{i,M}( \Delta t )\le 0 $ if $ \lambda_{i,o}\ge \lambda_{i,o1} $. Therefore, for any finite time $ t $, if $ \lambda_{i,o}\ge \lambda_{i,o1} $, then $ \dot{V}_{i,M}(t)\le 0 $ holds, which guarantees $ |\tilde{\rho}_{if}(t)|\le |\tilde{\rho}_{if}(0)|\le L_{0} $. 
	
	Next, with $-|\tilde{\rho}_{io}|\le -|\tilde{\rho}_{if}| $ and $-|(\beta_{1}\tilde{\rho}_{io}^{[r_{1}]}+\beta_{2}\tilde{\rho}_{io}^{[r_{2}]})^{[r_{0}]}|\le- |(\beta_{1}\tilde{\rho}_{if}^{[r_{1}]}+\beta_{2}\tilde{\rho}_{if}^{[r_{2}]})^{[r_{0}]}| $,   \eqref{d22} is rewritten as 
	\begin{align}\label{22}
	\dot{V}_{i,M}\le&-\gamma_{i,o}\lambda_{iv}|\tilde{\rho}_{io}||(\beta_{1}\tilde{\rho}_{io}^{[r_{1}]}+\beta_{2}\tilde{\rho}_{io}^{[r_{2}]})^{[r_{0}]}|\nm\\
	\le&-\frac{\gamma_{i,o}\lambda_{iv}}{2}|\tilde{\rho}_{io}||(\beta_{1}\tilde{\rho}_{io}^{[r_{1}]}+\beta_{2}\tilde{\rho}_{io}^{[r_{2}]})^{[r_{0}]}|\nm\\&-\frac{\gamma_{i,o}\lambda_{iv}}{2}|\tilde{\rho}_{if}||(\beta_{1}\tilde{\rho}_{if}^{[r_{1}]}+\beta_{2}\tilde{\rho}_{if}^{[r_{2}]})^{[r_{0}]}|,\nm\\ \le& -\eta_{iM1}{V}_{i,M}^{\frac{r_{1}r_{0}+1}{2}}-\eta_{iM2}{V}_{i,M}^{\frac{r_{2}r_{0}+1}{2}},
	\end{align}
	{where\\ {\small $ \eta_{iM1}:=\min\left\{\frac{\gamma_{i,o}\lambda_{iv}\beta_{1}}{2^{\frac{3-r_{1}r_{0}}{2}}}\left(\frac{2}{\gamma_{i,o}}\right)^{\frac{2}{r_{1}r_{0}+1}}, \frac{\gamma_{i,o}\lambda_{iv}\beta_{1}}{2^{\frac{3-r_{1}r_{0}}{2}}}\left(\frac{2}{\gamma_{f}}\right)^{\frac{2}{r_{1}r_{0}+1}} \right\} $, $ \eta_{iM2}:=\min\left\{\frac{\gamma_{i,o}\lambda_{iv}\beta_{2}}{2}\left(\frac{2}{\gamma_{i,o}}\right)^{\frac{2}{r_{2}r_{0}+1}}, \frac{\gamma_{i,o}\lambda_{iv}\beta_{2}}{2}\left(\frac{2}{\gamma_{f}}\right)^{\frac{2}{r_{2}r_{0}+1}} \right\}$,} and the following lemma is applied.
		\begin{lemma}\label{ld2}\cite{zuo2015nonsingular}
			Let $ \xi_{1},\ \xi_{2},\ \dots,\ \xi_{N}\ge 0 $, $ k_{1}>1 $, and $ 0<k_{2}\le1 $. Then 
			\[   \sum_{i=1}^{N}\xi_{i}^{k_{1}}\ge N^{1-k_{1}} \left(\sum_{i=1}^{N}\xi_{i}\right)^{k_{1}},\
			\sum_{i=1}^{N}\xi_{i}^{k_{2}}\ge \left(\sum_{i=1}^{N}\xi_{i}\right)^{k_{2}}.  \]
		\end{lemma}}
	 Then from Lemma~\ref{ld1}, for any {$(\tilde{\rho}_{io}(0),\tilde{\rho}_{if}(0))\in\Omega_{of}\times\Omega_{of}$ with $|\tilde{\rho}_{io}(0)|\ge |\tilde{\rho}_{if}(0)| > \phi_{s}$,} there exists a settling time 
	$$ T_{i,1}=\frac{2}{\eta_{iM1}(r_{1}r_{0}-1)}+\frac{2}{\eta_{iM2}(1-r_{2}r_{0})} $$
	such that $ \|x_{i1}-x_{i}^{o}\|\ge \sqrt{d^{2}-2\phi_{s}} $ and $ \|x_{i1}-\hat{x}_{i1}\|\ge \sqrt{d_{i0}^{2}-2\phi_{s}} $ for all $t\geq T_{i,1}$, {where $ \phi_{s} \in \mathbb{R}_{>0} $,\\$ \Omega_{of}:= (-\infty,-\phi_{s})\cup(\phi_{s},+\infty).$}
	\\ 
	
	{\textbf{(b)}} If $ \phi_{s}<|\tilde{\rho}_{io}|< |\tilde{\rho}_{if}|\le L_{0} {<+\infty}$, then $$\phi_{s}^{*}<|(\beta_{1}\tilde{\rho}_{io}^{[r_{1}]}+\beta_{2}\tilde{\rho}_{io}^{[r_{2}]})^{[r_{0}]}|< |(\beta_{1}\tilde{\rho}_{if}^{[r_{1}]}+\beta_{2}\tilde{\rho}_{if}^{[r_{2}]})^{[r_{0}]}|\le L_{0}^{*},$$
	where $ L_{0}^{*}:= (\beta_{1}L_{0}^{r_{1}}+\beta_{2}L_{0}^{r_{2}})^{r_{0}} $. Following the similar reasoning as the proof in \textbf{(a)}, we can verify that
	\begin{align*}
	\dot{V}_{i,M}(0)\le -\gamma_{i,o}\lambda_{iv}|\tilde{\rho}_{io}(0)||(\beta_{1}\tilde{\rho}_{io}^{[r_{1}]}(0)+\beta_{2}\tilde{\rho}_{io}^{[r_{2}]}(0))^{[r_{0}]}|\le 0,
	\end{align*} 
	where $ \lambda_{i,o} $ is designed to fulfill
	\begin{align}\label{d24}
	\nm \lambda_{i,o} \ge \lambda_{i,o2}:=& \frac{\gamma_{i,o}\lambda_{f}L_{iof}L_{0}^{*}}{\gamma_{i,\varepsilon}\phi_{s}^{*}}+\frac{\gamma_{i,o}L_{io}L_{i\hat{x}}}{\phi_{s}^{*}\gamma_{i,\varepsilon}}\\&+ \frac{\gamma_{f}L_{if}{(L_{i\hat{x}}+L_{i}^{o})}L_{0}}{\phi_{s}^{*}\phi_{s}\gamma_{i,\varepsilon}} +\frac{\gamma_{i,o}{\lambda_{iv}}}{\gamma_{i,\varepsilon}}.
	\end{align}
	Therefore, for any finite time $ t $, if $ \lambda_{i,o}\ge \lambda_{i,o2} $, then $ \dot{V}_{i,M}(t)\le 0 $ holds, which gives $ |\tilde{\rho}_{if}(t)|\le |\tilde{\rho}_{if}(0)|\le L_{0} $. 
	
	Moreover, using $-{|\tilde{\rho}_{io}(t)|}\le- \frac{\phi_{s}}{ L_{0}}{|\tilde{\rho}_{if}(t)|} $ and $ - |(\beta_{1}\tilde{\rho}_{io}^{[r_{1}]}+\beta_{2}\tilde{\rho}_{io}^{[r_{2}]})^{[r_{0}]}|\le -(\frac{\phi_{s}}{ L_{0}})^{r_{1}r_{0}} |(\beta_{1}\tilde{\rho}_{if}^{[r_{1}]}+\beta_{2}\tilde{\rho}_{if}^{[r_{2}]})^{[r_{0}]}| $,  we obtain 
	\begin{align}\label{d25}
	\dot{V}_{i,M}\le&-\eta_{iM3}{V}_{i,M}^{\frac{r_{1}r_{0}+1}{2}}-\eta_{iM4}{V}_{i,M}^{\frac{r_{2}r_{0}+1}{2}},
	\end{align}
	where
	$ \eta_{iM3}:=\\ \min \left\{\frac{\gamma_{i,o}\lambda_{iv}\beta_{1}}{2^{\frac{3-r_{1}r_{0}}{2}}}\left(\frac{2}{\gamma_{i,o}}\right)^{\frac{2}{r_{1}r_{0}+1}}, \frac{\gamma_{i,o}\lambda_{iv}\beta_{1}\phi_{s}\phi_{s}^{*}}{L_{0}L_{0}^{*}2^{\frac{3-r_{1}r_{0}}{2}}}\left(\frac{2}{\gamma_{f}}\right)^{\frac{2}{r_{1}r_{0}+1}} \right\} $, $ \eta_{iM4}:=\\ \min \left\{\frac{\gamma_{i,o}\lambda_{iv}\beta_{2}}{2}\left(\frac{2}{\gamma_{i,o}}\right)^{\frac{2}{r_{2}r_{0}+1}}, \frac{\gamma_{i,o}\lambda_{iv}\beta_{2}\phi_{s}\phi_{s}^{*}}{2L_{0}L_{0}^{*}}\left(\frac{2}{\gamma_{f}}\right)^{\frac{2}{r_{2}r_{0}+1}} \right\}$. 
	Then from  Lemma \ref{ld1}, for any {$(\tilde{\rho}_{io}(0),\tilde{\rho}_{if}(0))\in\Omega_{of}\times\Omega_{of}$ with $ \phi_{s}<|\tilde{\rho}_{io}(0)|< |\tilde{\rho}_{if}(0)|\le L_{0} $}, there exists a settling time
	$$ {T_{i,2}}:=\frac{2}{\eta_{iM3}(r_{1}r_{0}-1)}+\frac{2}{\eta_{iM4}(1-r_{2}r_{0})} $$ such that $ \|x_{i1}-x_{i}^{o}\|\ge \sqrt{d^{2}-2\phi_{s}} $ and $ \|x_{i1}-\hat{x}_{i1}\|\ge \sqrt{d_{i0}^{2}-2\phi_{s}} $ for all {$t\geq T_{i,2}$, where $\phi_{s}<L_{0} \in \mathbb{R}_{>0}$.}
	
	When the {local minima} is reached, i.e., $ \|\dot{x}_{io}-\dot{x}_{if}\|= 0 $, the same conclusions in \eqref{22} and \eqref{d25} are obtained by designing
	\begin{align*}
	\lambda_{i,o}\ge \lambda_{i,o3}:=& \frac{\gamma_{i,o}\lambda_{f}L_{iof}}{\gamma_{i,\varepsilon}}+ \frac{\gamma_{i,o}L_{io}(L_{i\hat{x}}+\delta_{d})}{\phi_{s}^{*}\gamma_{i,\varepsilon}}\nm\\& + \frac{\gamma_{f}L_{if}(L_{i\hat{x}}+L_{i}^{o}+\delta_{d})}{\phi_{s}^{*}\gamma_{i,\varepsilon}} +\frac{\gamma_{i,o}\gamma_{i,v}}{\gamma_{i,\varepsilon}},
	\\
	\nm \lambda_{i,o}\ge \lambda_{i,o4}:=& \frac{\gamma_{i,o}\lambda_{f}L_{iof}L_{0}^{*}}{\gamma_{i,\varepsilon}\phi_{s}^{*}}+\frac{\gamma_{i,o}L_{io}(L_{i\hat{x}}+
		\delta_{d})}{\phi_{s}^{*}\gamma_{i,\varepsilon}}\\&+ \frac{\gamma_{f}L_{if}(L_{i\hat{x}}+L_{i}^{o}+\delta_{d})L_{0}}{\phi_{s}^{*}\phi_{s}\gamma_{i,\varepsilon}} +\frac{\gamma_{i,o}\gamma_{i,v}}{\gamma_{i,\varepsilon}}.
	\end{align*}
	
	In conclusion, for any $(\tilde{\rho}_{io}(0),\tilde{\rho}_{if}(0))\in\bR\times\bR$ and {$ \lambda_{i,o}\ge \max\{\lambda_{i,o1},\lambda_{i,o2},\lambda_{i,o3},\lambda_{i,o4}\} $}, if we utilize the merged desired velocity in \eqref{v} for $\tilde{\rho}_{io}$ and $\tilde{\rho}_{if}$, then there exists a settling time $$T_{i}:=\max\{T_{i,o},T_{i,f},T_{i,1},T_{i,2}\}$$ such that $ \|x_{i1}-x_{i}^{o}\|\ge \sqrt{d^{2}-2\phi_{s}} $ and $ \|x_{i1}-\hat{x}_{i1}\|\ge \sqrt{d_{i0}^{2}-2\phi_{s}} $ for all $t\geq T_{i}$.
	
	{In addition, when $ \|\dot{x}_{io}+\dot{x}_{if}\|= 0 $, the proof is similar to the case $ \|\dot{x}_{io}-\dot{x}_{if}\|\neq 0 $. The only difference is that $ L_{i\hat{x}}$ is redefined as $ L_{i\hat{x}}:= K_{1}\|\bar{\eta}_{i1}(0)\|+K_{2}\|\bar{\eta}_{i2}(0)\|+\sqrt{3}K_{3} +\delta_{d}$.}

	\subsection{Proof of Theorem~\ref{the1}}\label{Ap:the}
	
	The proof of this theorem contains two parts. In the first part, we show that $ \|{\mathcal{S}}_{i}(\cdot)\|$, $\|\tilde{W}_{i}\|_F $ and $\|\tilde{\delta}_{i}\| $ are bounded, and based on it, we prove the fixed-time convergence of the tracking errors ${\tilde{x}}_{i1}$, $\tilde{x}_{i2}$ in the second part.
	
	\textbf{(1)} \emph{Boundedness of $\mathcal{S}_{i}$, $\tilde{W}_{i}$ and $\tilde{\delta}_{i}$.}
	
	Consider the Lyapunov candidate as follows:
	\begin{align*} V(\mathcal{S},\tilde{W},\tilde{\delta}):&=V_{S}(\mathcal{S})+V_{P1}(\tilde{W},\tilde{\delta}), \ \ {\rm with}\\
	V_{S}(\mathcal{S}):&=\sum\limits_{i=1}^{n}\frac{1}{2\varrho}\mathcal{S}_{i}^{\top}\mathcal{S}_{i}, \\
	V_{P1}(\tilde{W},\tilde{\delta}): &= \frac{1}{2} \sum\limits_{i=1}^{n}{\rm Tr}[\tilde{W}_{i}^{\top}\Gamma_{i}^{-1}\tilde{W}_{i}]+\sum\limits_{i=1}^{n}\frac{1}{2\gamma_{3i}}\tilde{\delta}_{i}^2,
	\end{align*}
	where $ \varrho>0$ is a design parameter, ${\mathcal{S}}:=[{\mathcal{S}}_{1}^{\top},\dots,{\mathcal{S}}_{n}^{\top}]^{\top} $, $ \tilde{W}:={\rm blockdiag}\{\tilde{W}_{1},\dots,\tilde{W}_{n}\} $ with $ \tilde{W}_{i}:={W}_{i}^{*}-\hat{W}_{i} $, $ \tilde{\delta}:=[\tilde{\delta}_{1},\dots,\tilde{\delta}_{n}]^{\top} $, $ {\rm Tr}(A) $ denotes the trace of a square matrix $ A $, $ \Gamma_{i}{\in\bR^{h_{i}\times h_{i}}} $ is a positive definite constant gain matrix.
	
	With the differentiation of \eqref{S}, we calculate the time derivative of $ {V}_{S}(\mathcal{S}) $ as  
	\begin{align} 
	\dot{V}_{S}
	=& \sum\limits_{i=1}^{n}{\mathcal{S}}_{i}^{\top} \left[ \dot{x}_{i2}-\ddot{x}_{id}+c_{1}\dot{\tilde{x}}_{i1}+c_{2}\dot{\alpha}_{i}(\dot{\tilde{x}}_{i1})\right ]
	\nm\\
	= & \sum\limits_{i=1}^{n}{\mathcal{S}}_{i}^{\top} \left[ f_{i}(\bar{x}_{i})+u_{i}+d_{i}+\chi_{i}\right ],
	\label{eq:dotVs} 
	\end{align}
	with $u_{i}$ in \eqref{u} and  $\chi_{i}(z_{1})=-\ddot{x}_{id}+c_{1}\dot{\tilde{x}}_{i1}+c_{2}\dot{\alpha}_{i}(\dot{\tilde{x}}_{i1})$. 
	Let $ z_{i}=:[z_{1i}^{\top},\dots,z_{1i}^{\top}]^{\top} $, with $ z_{1i}:=[\ddot{x}_{id}^{\top},\dot{\tilde{x}}_{i1}^{\top},\dot{\alpha}_{i}(\dot{\tilde{x}}_{i1})^{\top}]^{\top} $. The uncertainty term $ F_{i}(z_{i}):=f_{i}(\bar{x}_{i})+\chi_{i}(z_{1i}) $ in \eqref{32} is estimated by the RBFNN as $ F_{i}(z_{i}):={W_{i}^{\ast}(t)}^{\top}\phi_{i}(z_{i})+\varepsilon_{i},\ \forall t\in\bR_{\ge 0},z_{i}\in\Omega_{zi} $ in \eqref{33}. As a result, \eqref{eq:dotVs} becomes
	\begin{align}
	\dot{V}_{S} = &\sum\limits_{i=1}^{n}{\mathcal{S}}_{i}^{\top} [ W_{i}^{\ast\top}\phi_{i}(z_{i})+\varepsilon_{i}+d_{i} +u_{i}^{S}+u_{i}^N+u_{i}^C]
	\nm\\
	= & \sum\limits_{i=1}^{n}{\mathcal{S}}_{i}^{\top} [ W_{i}^{\ast\top}\phi_{i}(z_{i})-\hat{W}_{i}^{\top}\phi_{i}(z_{i})+\varepsilon_{i}+d_{i} -\hat{\delta}_{i}\text{sgn}({\mathcal{S}}_{i})\nm\\
	& - k_{1i}{\mathcal{S}}_{i}^{[\gamma_{1}]}-k_{2i}{\mathcal{S}}_{i}^{[\gamma_{2}]}]\nonumber\\
	\label{eq:dotVs2}
	%\dot{V}_{S} = &\sum\limits_{i=1}^{n}{\mathcal{S}}_{i}^{\top} [ W_{i}^{\ast\top}\phi_{i}(z_{i})+\varepsilon_{i}+d_{i} +u_{i}^{S}+u_{i}^N+u_{i}^C]
	%\nm\\
	%= & \sum\limits_{i=1}^{n}{\mathcal{S}}_{i}^{\top} [ W_{i}^{\ast\top}\phi_{i}(z_{i})-\hat{W}_{i}^{\top}\phi_{i}(z_{i})+\varepsilon_{i}+d_{i} -\hat{\delta}_{i}\text{sgn}({\mathcal{S}}_{i}) \nm\\
	%& - k_{1i}{\mathcal{S}}_{i}^{[\gamma_{1}]}-k_{2i}{\mathcal{S}}_{i}^{[\gamma_{2}]}]
	%\nm\\
	\leq & \sum\limits_{i=1}^{n}{\mathcal{S}}_{i}^{\top} [-k_{1i}{\mathcal{S}}_{i}^{[\gamma_{1}]}-k_{2i}{\mathcal{S}}_{i}^{[\gamma_{2}]}] + \sum\limits_{i=1}^{n}(\delta_{i}-\hat{\delta}_{i})\|{\mathcal{S}}_{i}\|_{1}
	\nm\\ &+ \sum\limits_{i=1}^{n} {\mathcal{S}}_{i}^{\top}(W_{i}^{\ast}-\hat{W}_{i})^{\top}\phi_{i}(z_{i}),
	\end{align}
	where the assumption $\|\varepsilon_{i}+d_{i}\|\leq \delta_{i}$ is used. 
	
	Then we further compute the derivative of $ V_{P1}(\tilde{W},\tilde{\delta}) $. From \eqref{ad1} and \eqref{ad2}, it follows that
	\begin{align}\label{54}
	\dot{V}_{P1} =& -\sum\limits_{i=1}^{n}\text{Tr}[\tilde{W}_{i}^{\top}\phi_{i}(z_{i}){\mathcal{S}}_{i}^{\top}]-\sum\limits_{i=1}^{n}\|\tilde{\delta}_{i}\|_{1}{\mathcal{S}}_{i}.
	\end{align}
	with $\tilde{W}_{i}:=W_{i}^{\ast}-\hat{W}_{i}$, and $\tilde{\delta}_{i}:=\delta_{i}-\hat{\delta}_{i}$. Note that $\text{Tr}[\tilde{W}_{i}^{\top}\phi_{i}(z_{i}){\mathcal{S}}_{i}^{\top}] = {\mathcal{S}}_{i}^{\top} \tilde{W}_{i}^{\top} \phi_{i}(z_{i})$. Therefore, 
	\begin{align*}
	\dot{V}(\mathcal{S},\tilde{W},\tilde{\delta}) &=\dot{V}_{S}(\mathcal{S})+\dot{V}_{P1}(\tilde{W},\tilde{\delta}) \\
	&\leq -\sum\limits_{i=1}^{n}k_{1i}{\mathcal{S}}_{i}^{\top}{\mathcal{S}}_{i}^{[\gamma_{1}]}-\sum\limits_{i=1}^{n}k_{2i}{\mathcal{S}}_{i}^{\top}{\mathcal{S}}_{i}^{[\gamma_{2}]}\leq 0,
	\end{align*}  
	which implies that for any initial conditions $ {\mathcal{S}}_{i}(0),\hat{W}_{i}(0),\hat{\delta}_{i}(0)$, there exist positive constants $\varsigma_{si} $, $ \varsigma_{0i} $, $ \varsigma_{1i} $, which depend on the values of $ {\mathcal{S}}_{i}(0),\hat{W}_{i}(0),\hat{\delta}_{i}(0)$, such that $ \|{\mathcal{S}}_{i}(\cdot)\|\le \varsigma_{si} $, $\|\tilde{W}_{i}\|_F \leq \varsigma_{0i}$ and $\|\tilde{\delta}_{i}\| \leq \varsigma_{1i}$. \\
	
	\textbf{(2)} \emph{Fixed-time convergence of the tracking errors ${\tilde{x}}_{i1}$, $\tilde{x}_{i2}$.}
	
	Since $0<\phi_{ki}(z_{i})\le 1$ holds for all the neurons $k = 1, 2, \cdots, h_{i}$, we have 
	$\|\phi_{i}(z_{i})\|\leq \sqrt{h_{i}}$ as
	$ \phi_{i}(z_{i})=[\phi_{1i}(z_{i}),\phi_{2i}(z_{i}),\ldots,\phi_{h_{i}i}(z_{i})]^{\top}$. Then using the property of the Frobenius norm, we obtain   $$\|\tilde{W}_{i}^{\top}\phi_{i}(z_{i})\|_F \leq \|\tilde{W}_{i}\|_F\|\phi_{i}(z_{i})\| = \sqrt{h_{i}}\varsigma_{0i}. 
	$$
	Denote $ \varsigma_{2i}:=\sqrt{h_{i}}\varsigma_{0i} $ and $\varsigma_{3i}=\text{max}\{\varsigma_{1i}+\varsigma_{2i}\}$. With the property $\|{\mathcal{S}}_{i}\|_{1}\le \sqrt{3} \|{\mathcal{S}}_{i}\|,$ it follows from \eqref{eq:dotVs2} that
	\begin{align}
	\dot{V}_{S} &\leq -\sum\limits_{i=1}^{n}k_{1i}{\mathcal{S}}_{i}^{\top}{\mathcal{S}}_{i}^{[\gamma_{1}]}-\sum\limits_{i=1}^{n}k_{2i}{\mathcal{S}}_{i}^{\top}{\mathcal{S}}_{i}^{[\gamma_{2}]} +\sum\limits_{i=1}^{n}\varsigma_{3i}\|{\mathcal{S}}_{i}\|_{1} 
	\nm\\
	&\le -\gamma_{s1} V_{S}^{\frac{1+\gamma_{1}}{2}}-\gamma_{s2}V_{S}^{\frac{1+\gamma_{2}}{2}}+\gamma_{s3}V_{S}^{\frac{1}{2}}, \label{eq:dotVs3}
	\end{align}
	where $ \gamma_{s1}:=k_{1i}(2\varrho)^{\frac{1+\gamma_{1}}{2}} $, $ \gamma_{s2}:=k_{2i}(2\varrho)^{\frac{1+\gamma_{2}}{2}} $, $ \gamma_{s3}:=\sqrt{3}(2\varrho)^{2}\max_{i}\{\varsigma_{3i}\}$.

	To eliminate the third positive term in \eqref{eq:dotVs3}, we write $ \gamma_{s1}:=\gamma_{s1} ^{a}+\gamma_{s1} ^{b} $ with $ \gamma_{s1} ^{a}$, $\gamma_{s1} ^{b} \in \mathbb{R}_{>0}$ such that
	$$
	\|{\mathcal{S}}_{i}\|\ge \delta_{s3}:= \sqrt{2\varrho}(\gamma_{s3}/\gamma_{s1} ^{b})^\frac{1}{\gamma_{1}}.$$
	Then we obtain
	$ \dot{V}_{S}\le -\gamma_{s1} ^{a} V_{S}^{\frac{1+\gamma_{1}}{2}}-\gamma_{s2}V_{S}^{\frac{1+\gamma_{2}}{2}}. $
	Alternatively, we can write $ \gamma_{s2}=\gamma_{s2} ^{a}+\gamma_{s2} ^{b} $ with $ \gamma_{s2} ^{a}$, $\gamma_{s2} ^{b} \in \mathbb{R}_{>0}$, and
	$$ \|{\mathcal{S}}_{i}\|\ge \delta_{s4}:= \sqrt{2\varrho}(\gamma_{s3}/\gamma_{s2} ^{b})^\frac{1}{\gamma_{2}}, $$
	which leads to
	$ \dot{V}_{S}\le -\gamma_{s1} V_{S}^{\frac{1+\gamma_{1}}{2}}-\gamma_{s2} ^{a}V_{S}^{\frac{1+\gamma_{2}}{2}}. 
	$
	Combine the two ways, we define $ \delta_{s}:=\min\{\delta_{s3},\delta_{s4}\} $. Note that $ \delta_{s} $ can be made arbitrary small by making $ k_{1i} (t)$ and $ k_{2i} (t)$ sufficiently large. Therefore, according to Lemma~\ref{ld1}, there exists  
	$ T_{S1}\le\max\{T_{\max}^{a},T_{\max}^{b}\}, $ where
	\begin{align*}
	T_{\max}^{a}:=\frac{1}{\gamma_{s1} ^{a}(\frac{1+\gamma_{1}}{2}-1)}+\frac{1}{\gamma_{s2}(1-\frac{1+\gamma_{2}}{2})}, \\
	T_{\max}^{b}:=\frac{1}{\gamma_{s1}(\frac{1+\gamma_{1}}{2}-1)}+\frac{1}{\gamma_{s2} ^{a}(1-\frac{1+\gamma_{2}}{2})},
	\end{align*}
	such that $ \|{\mathcal{S}}_{i}\|\le \delta_{s} $ for all $ t\ge T_{S1}>0 $.
	This completes the proof.

\subsection{Proof of Lemma~\ref{l7}}\label{Ap:l7}

For any $ \delta_{s} \in \mathbb{R}_{>0}$, if $ \|{\mathcal{S}}_{i}(\cdot)\|\le\delta_{s} $, then $ |{\mathcal{S}}_{i_{j}}(\cdot)|\le\delta_{s} $ for $ j=1,2,3 $. Three cases are discussed based on the definition of $ {\mathcal{S}}_{i_{j}}(\cdot) $ in \eqref{S}.

{\bf Case A}: $ \sigma_{1i_{j}}(\cdot) =0$ for all $ j=1,2,3 $. According to Lemma~6, there exists a fixed $ T_{s_{0}} >0$ such that $ \lim_{t\to T_{s_{0}}} \tilde{x}_{i1_{j}}(t)=0$, $ \lim_{t\to T_{s_{0}}} \tilde{x}_{i2_{j}}(t)=0$. 

{\bf Case B}: $ \sigma_{1i_{j}} (\cdot)\neq0,\ |\tilde{x}_{i1_{j}}| \leq \phi_{s} $ for some $ j $. It follows from the definition of $ {\mathcal{S}}_{i_{j}}(\cdot)$ that
$$ |\tilde{x}_{i2_{j}}+c_{1}\tilde{x}_{i1_{j}}+c_{2}(\wp_{1}\tilde{x}_{i1_{j}}+\wp_{2}\tilde{x}_{i1_{j}}^{[2]})|\le \delta_{s}/\varrho. $$
Then from \eqref{8}, \eqref{9}, and $ |\tilde{x}_{i1_{j}}| \leq \phi_{s} $, we obtain
\begin{equation}
	|\tilde{x}_{i2_{j}}|\le \delta_{s}/\varrho+ c_{1}\phi_{s}+c_{2}(\beta_{1}\phi_{s}^{r_{1}}+\beta_{2}\phi_{s}^{r_{2}})^{r_{0}}:=\delta_{s_{0}}.
\end{equation}

{\bf Case C}: $ \sigma_{1i_{j}} (\cdot)\neq0,\ |\tilde{x}_{i1_{j}}| > \phi_{s} $ for some $ j $. We have $$| \tilde{x}_{i2_{j}}+c_{1}\tilde{x}_{i1_{j}}+c_{2}\mathcal{B}^{[r_{0}]} |\le \delta_{s}/\varrho, $$
where $\mathcal{B}: = \beta_{1}\tilde{x}_{i1_{j}}^{[r_{1}]}+\beta_{2}\tilde{x}_{i1_{j}}^{[r_{2}]}$.  
If $0\le \tilde{x}_{i2_{j}}+c_{1}\tilde{x}_{i1_{j}}+c_{2}\mathcal{B}^{[r_{0}]} \le \delta_{s}/\varrho$,
we consider the Lyapunov function 
$ V_{S1}(\tilde{x}_{i1_{j}}):=\frac{1}{2}\tilde{x}_{i1_{j}}^{2}, 
$ whose time derivative is given as
\begin{align*}
	\dot{V}_{S1}=\tilde{x}_{i1_{j}}\tilde{x}_{i2_{j}} 
	&\le \tilde{x}_{i1_{j}}\left[-c_{1}\tilde{x}_{i1_{j}}-c_{2}\mathcal{B}^{[r_{0}]}+\delta_{s}/\varrho\right]
	\\
	&\le \tilde{x}_{i1_{j}}\left[-(c_{1}-{\delta_{s}}/{\varrho\tilde{x}_{i1_{j}}})\tilde{x}_{i1_{j}}-c_{2}\mathcal{B}^{[r_{0}]}\right]   
	\\
	& \le 
	\tilde{x}_{i1_{j}}\left[ - c_{1}\tilde{x}_{i1_{j}}+(c_{2}-{\delta_{s}}/{(\varrho\mathcal{B}^{[r_{0}]})})
	\mathcal{B}^{[r_{0}]}\right]. 
\end{align*}
Consider the parameters $\tilde{c}_{1}$,  $\tilde{c}_{2}$, $\tilde{\beta}_{1}$, $\tilde{\beta}_{2}$ defined in this lemma, and we denote 
$c_{1}^a: =  c_{1} - \tilde{c}_{1} $, 
$ c_{2}^a =  c_{2} - \tilde{c}_{2} $, 
$ \beta_{1}^a = \beta_{1} - \tilde{\beta}_{1} $, 
$ \beta_{2}^a = \beta_{2} - \tilde{\beta}_{2} $, 
which are positive scalars. Using the three inequalities in \eqref{eq:tildecbeta}, we obtain
\begin{subequations}\label{fixed3}
	\begin{align}
		\dot{V}_{S1}\le -c_{1}^a\tilde{x}_{i1_{j}}^{2}-c_{2}(\beta_{1}\tilde{x}_{i1_{j}}^{\frac{r_{1}+1}{r_{0}}}+\beta_{2}\tilde{x}_{i1_{j}}^{\frac{r_{2}+1}{r_{0}}})^{r_{0}},
		\\ 
		\dot{V}_{S1}\le - c_{1}\tilde{x}_{i1_{j}}^{2}-c_{2}^a(\beta_{1}^a\tilde{x}_{i1_{j}}^{\frac{r_{1}+1}{r_{0}}}+\beta_{2}\tilde{x}_{i1_{j}}^{\frac{r_{2}+1}{r_{0}}})^{r_{0}},
		\\
		\dot{V}_{S1}\le- c_{1}\tilde{x}_{i1_{j}}^{2}-c_{2}^a(\beta_{1}\tilde{x}_{i1_{j}}^{\frac{r_{1}+1}{r_{0}}}+\beta_{2}^a\tilde{x}_{i1_{j}}^{\frac{r_{2}+1}{r_{0}}})^{r_{0}},
	\end{align}
\end{subequations}
respectively. If one of the above inequalities holds, we then conclude from Remark~\ref{remark2} that: for any $|\tilde{x}_{i1_{j}}(0)|>0 $, there exists a fixed time $ T_{s_{1}} >0$ such that 
$$ |\tilde{x}_{i1_{j}}(t)|\le\min\left\{
\frac{\delta_{s}}{\varrho \tilde{c}_{1}},
\frac{(\delta_{s}/\varrho \tilde{c}_{2})^{1/r_{0}r_{1}}}{(\tilde{\beta}_{1})^{1/r_{1}}},
\frac{(\delta_{s}/ \varrho \tilde{c}_{2})^{1/r_{0}r_{2}}}{(\tilde{\beta}_{2})^{1/r_{2}}}\right\} : = \delta_{s_{1}}
$$
for all $ t\ge T_{s_{1}} $. Furthermore, we get  \begin{equation}
	|\tilde{x}_{i2_{j}}|\le \delta_{s}/ \varrho+ c_{1}\delta_{s_{1}}+c_{2}(\beta_{1}\delta_{s_{1}}^{[r_{1}]}+\beta_{2}\delta_{s_{1}}^{[r_{2}]})^{[r_{0}]} :=\delta_{s_{2}}.
\end{equation} 

If $-\delta_{s}/\varrho \le \tilde{x}_{i2}+c_{1}\tilde{x}_{i1}+c_{2}(\beta_{1}\tilde{x}_{i1}^{[r_{1}]}+\beta_{2}\tilde{x}_{i1}^{[r_{2}]})^{[r_{0}]} \le 0 $, we can use the similar reasoning to consider a Lyapunov function $ V_{S1}(\tilde{x}_{i1_{j}}):=\frac{1}{2}\tilde{x}_{i1_{j}}^{2} $, and it can be verified that
$$ \dot{V}_{S1}\le -c_{1}\tilde{x}_{i1_{j}}^{2}-c_{2}(\beta_{1}\tilde{x}_{i1_{j}}^{\frac{r_{1}+1}{r_{0}}}+\beta_{2}\tilde{x}_{i1_{j}}^{\frac{r_{2}+1}{r_{0}}})^{r_{0}}. $$
Then from Remark~\ref{remark2}, it is obvious that for any $|\tilde{x}_{i1_{j}}(0)|>0 $, there exists a fixed time $ T_{s_{2}} >0$ such that 
$ \lim_{t\to T_{s_{2}} }\tilde{x}_{i1_{j}}(t)=0$ and $ \lim_{t\to T_{s_{2}} }\tilde{x}_{i2_{j}}(t)=0$
for all $ t\ge T_{s_{2}} $.
In conclusion, for any $|\tilde{x}_{i1_{j}}(0)|>0 $, there exists a fixed time $ T_{s_{3}}:=\max\{T_{s_{1}},T_{s_{2}}\} >0$ such that 
$ |\tilde{x}_{i1_{j}}(t)|\le\delta_{s_{1}}$, $ |\tilde{x}_{i2_{j}}|\le \delta_{s_{2}}$
for all $ t\ge T_{s_{3}} $.

Finally, from the analysis of the three cases, we conclude that for any $|\tilde{x}_{i1_{j}}(0)|,|\tilde{x}_{i2_{j}}(0)|>0 $, there exists a fixed time $ T_{s_{*}}:=\max\{T_{s_{0}},T_{s_{1}},T_{s_{2}}\} >0$ such that 
$ |\tilde{x}_{i1_{j}}(t)|\le\max \{\phi_{s},\delta_{s_{1}}\}$, $ |\tilde{x}_{i2_{j}}(t)|\le \max\{\delta_{s_{0}},\delta_{s_{2}}\}$
for all $ t\ge T_{s_{*}} $.	
	
	\bibliographystyle{IEEEtran}
	\bibliography{ref}

\end{document}